\documentclass[reqno,11pt]{amsart}
\usepackage{amsmath,amssymb,mathrsfs,amsthm,amsfonts}
\usepackage[inline]{enumitem}
\usepackage[usenames,dvipsnames]{xcolor}
\usepackage{hyperref}
\usepackage[percent]{overpic}
\usepackage{comment}
\usepackage{stmaryrd}
\usepackage{algorithm}
\usepackage{algorithmic}
\usepackage{pgfplots}
\usepackage{subcaption}
\usepackage{tikz,todonotes}
\usetikzlibrary{calc}
\usetikzlibrary{arrows.meta,backgrounds}
\usepackage{multirow,array,longtable,booktabs}
\usetikzlibrary{arrows}
\usepackage{bm}
\usepackage{graphicx}
\usepackage[T1]{fontenc}

\hypersetup{
	colorlinks=true, linkcolor=blue,
	citecolor=ForestGreen
}
\usepackage[paper=letterpaper,margin=1in]{geometry}
\DeclareMathAlphabet{\mathpzc}{OT1}{pzc}{m}{it}

\newtheorem{theorem}{Theorem}[section]

\newtheorem{lemma}{Lemma}[section]
\newtheorem{proposition}{Proposition}[section]

\newtheorem{definition}{Definition}[section]

\newtheorem{remark}{Remark}[section]
\theoremstyle{definition}
\newtheorem{example}{Example}[section]
\numberwithin{equation}{section}
\allowdisplaybreaks
\usepackage{acronym}

\newcommand{\be}{\begin{equation}}
\newcommand{\ee}{\end{equation}}

\makeatletter

\newcommand{\R}{\mathbb{R}} % real numbers
\newcommand{\e}{\varepsilon}

\newcommand{\T}{\mathbb{T}}

\newcommand{\nd}{\noindent}

\renewcommand{\hat}{\widehat}

\renewcommand{\bar}{\overline}

\usepackage{comment}

\newcounter{hyp}
\title[Long-time behavior of Hamilton-Jacobi equations in the Wasserstein space]{Long-time behavior of first-order Hamilton-Jacobi equations in the Wasserstein space}

\author[E.\ Bayraktar]{Erhan Bayraktar}
\address{E.\ Bayraktar,
	Department of Mathematics, University of Michigan,
	\newline\hphantom{\quad \ \ E. Bayraktar}
	530 Church St, Ann Arbor, MI 48109, USA}
\email{erhan@umich.edu}
\thanks{E. Bayraktar is partially supported by the NSF grant DMS-2507940, and in part by the Susan M. Smith Professorship.}

\author[A.\ Zitridis]{Antonios Zitridis}
\address{A.\ Zitridis,
	Department of Mathematics, University of Michigan,
	\newline\hphantom{\quad \ \ A. Zitridis}
	530 Church St, Ann Arbor, MI 48109, USA}
\email{zitridis@umich.edu}

\thanks{}

\begin{document}
\null
\vspace{-5mm}
    \maketitle
\vspace{-8mm}
\begin{abstract}
    We establish the long-time behavior of solutions $U$ to first-order Hamilton-Jacobi equations in the Wasserstein space $\mathcal{P}_2(\mathbb T^d)$. For a class of Hamiltonians with quadratic growth, including non-convex ones, we prove that there exist a unique ergodic constant $\lambda$ and a solution $V$ of the associated ergodic problem such that $U(t,m)-\lambda t$ converges to $V(m)$, uniformly on $\mathcal{P}_2(\mathbb T^d)$, as $t\to+\infty$. Finally, we provide applications to mean field control problems and zero-sum differential games.
    %This extends the classical large-time theory for Hamilton--Jacobi equations from finite-dimensional state spaces to the space of probability measures.
\end{abstract}

%\tableofcontents

\section{Introduction}
 
\nd
Let $d\in\mathbb{N}$ and $\T^d=\R^d/\mathbb{Z}^d$ be the $d$-dimensional torus. We consider a Hamiltonian $H\colon\T^d\times\R^d\times\mathcal{P}_2\to\R$ and a function $G\colon\mathcal{P}_2\to\R$, where $\mathcal{P}_2:=\mathcal{P}_2(\T^d)$ is the space of probability measures on $\T^d$ endowed with the $2$-Wasserstein distance ${\bf d}_2$. This paper is concerned with the long-time behavior of the solution $U\colon[0,+\infty)\times\mathcal{P}_2\to\R$ of the Hamilton--Jacobi--Bellman (HJ for short) equation
\be\label{HJB0}\tag{HJB}
\begin{cases}
    \partial_tU+\int_{\T^d}H\big(x,D_mU(t,m,x),m\big)\,m(dx)=0,& (t,m)\in(0,+\infty)\times\mathcal{P}_2,\\
    U(0,m)=G(m), &m\in \mathcal{P}_2,
\end{cases}
\ee
where the derivative $D_mU\colon[0,+\infty)\times\mathcal{P}_2\times\T^d\to\R^d$ is understood as the Wasserstein derivative (see \cite{carmona2018probabilistic}). Throughout the paper, we assume that $H$ has quadratic growth in $p$, but the map $p\mapsto H(x,p,m)$ is not necessarily convex.
 
\subsection{Hamilton--Jacobi equations in the Wasserstein space}
HJ equations posed on spaces of probability measures have attracted considerable attention in recent years due to their connections with several areas including mean field games, mean field control, optimal transport and filtering. For a sample of recent works, we refer to \cite{cardaliaguet2019master} for mean field games and mean field control, \cite{bertucci2024stochastic} for optimal transport, and \cite{martini2023kolmogorov,BayraktarEkrenHeZhang2025SC} for filtering.
\vspace{1mm}
 
\nd
In particular, when $H$ is convex in $p$, equations of the form \eqref{HJB0} are related to dynamic programming equations associated with mean field control problems. These problems study the collective behavior of a large population of interacting agents/particles who cooperate to minimize a common cost functional that might depend on the distribution of the population; see \cite{carmona2018probabilistic,cardaliaguet2019master} for surveys on the topic.
To make this connection explicit, assume for a moment that $p\mapsto H(x,p,m)$ is convex, and let $L(x,a,m):=\sup_{p\in\R^d}\{-a\cdot p-H(x,p,m)\}$, so that $H(x,p,m)=\sup_{a\in\R^d}\{-a\cdot p-L(x,a,m)\}$. For $T>0$, the function $W^T(t,m):=U(T-t,m)$ solves the backward equation $-\partial_tW^T+\int H(x,D_mW^T,m)\,dm=0$ in $(0,T)\times\mathcal{P}_2$ with terminal condition $W^T(T,\cdot)=G$, and it is, at least formally, the value function of the deterministic mean field control problem
\be\label{MFCvalue}
W^T(t,m):=\inf_{(\alpha,m)}\left\{\int_t^T\int_{\T^d}L\big(x,\alpha_s(x),m_s\big)\,m_s(dx)\,ds+G(m_T)\right\}.
\ee
Here, the infimum is taken over measurable $\alpha\colon[t,T]\times\T^d\to\R^d$ so that $(m_s)_{s\in[t,T]}$ solves the continuity equation $\partial_sm_s+\mathrm{div}(\alpha_sm_s)=0$ with $m_t=m$, and $\int_t^T\int |\alpha_s(x)|^2 m_s(dx)ds<\infty$ \footnote{Here we used the feedback control representation of \eqref{MFCvalue}. We refer to \cite[Remarks 1.4, 1.5]{ding2026quantitative} for a more detailed treatment of this minimization problem with different types of controls.}. We refer to \cite{gangbo2021finite} for the connection between \eqref{MFCvalue} and the HJ equation satisfied by $W^T$. In particular, $U(T,m)=W^T(0,m)$ is the optimal cost of controlling a population with initial distribution $m$ over the horizon $[0,T]$.
\vspace{1mm}
 
\nd
When $H$ is not convex in $p$, a representation such as \eqref{MFCvalue} is not available. Non-convex Hamiltonians arise, for instance, in zero-sum stochastic differential games. These games were originally formulated rigorously by
Elliott-Kalton \cite{elliott1972existence}, and their connection with viscosity solutions of HJ equations goes back to Evans-Souganidis \cite{evans1984differential}; see also Bardi-Capuzzo-Dolcetta \cite{bardi1997optimal} for the classical theory of deterministic control and differential games. For a weak dynamic programming principle and viscosity characterization for stochastic differential games with unbounded controls, see Bayraktar-Yao \cite{bayraktar2013weak}. For zero-sum games in infinite dimensions (games of McKean--Vlasov type), we refer to
Cosso-Pham \cite{cosso2019zero} (stochastic dynamics) and to Ba\c{s}ar-Moon \cite{moon2021zero} (deterministic games on the
Wasserstein space). We refer to section \ref{ex} for a more precise connection between \eqref{HJB0} and zero-sum differential McKean-Vlasov games.
%of McKean--Vlasov type, whose value functions solve Bellman--Isaacs equations in the Wasserstein space \cite{cosso2019zero}.
\vspace{1mm}
 
\nd
As in finite dimensions, HJ equations in the Wasserstein space do not have classical solutions in general\footnote{In our setup, this is already the case when $H$ does not depend on $m$ and $G$ is linear in $m$, since \eqref{HJB0} then reduces, at least formally, to a finite-dimensional HJ equation (see \eqref{HJfd} below), whose solutions typically develop singularities in finite time.}. Therefore, a large number of recent works on the topic have been devoted to developing an appropriate notion of viscosity sub/super-solution and to establishing a comparison principle. We refer to \cite{conforti2023hamilton,daudin2024comparison,bayraktar2025comparison,bayraktarqiu2025viscosity,cheung2025viscosity,daudin2025well,zhou2024viscosity,soner2024viscosity}, among others. Unlike in finite dimensions \cite{crandall1992user}, there are several reasonable definitions of viscosity solutions for equations posed on spaces of probability measures. Quantitative particle approximation for second-order equations, including non-convex Hamiltonians and common noise, is studied in \cite{bayraktar2025particle}. Most relevant to \eqref{HJB0} are the works \cite{gangbo2019differentiability,gangbo2021finite}, which consider viscosity solutions in the ``lifted'' sense of P.-L.~Lions: a function on $\mathcal{P}_2$ is identified with its lift to the Hilbert space $L^2(\Omega;\R^d)$ of random variables, and \eqref{HJB0} is replaced by an HJ equation on this Hilbert space. In this paper, we adopt this notion of viscosity solutions and we refer to them as $L$-viscosity solutions. We defer their precise definition to Section~\ref{viscositysol}, Definition~\ref{visc.}, where we also show that an $L$-viscosity solution of \eqref{HJB0} is obtained as the limit of solutions of finite-dimensional (particle) approximations, in the spirit of \cite{gangbo2021finite} (see also \cite{daudin2025error} for error estimates for such approximations). The uniqueness of $L$-viscosity solutions follows from a comparison principle (Theorem~\ref{comp}).
\vspace{2mm}

\nd
The goal of this paper is to describe the behavior of the $L$-viscosity solution $U(t,m)$ of \eqref{HJB0}, as $t\to+\infty$. More specifically, we show that $U(t,m)/t$ converges to a constant $\lambda$, as $t\to+\infty$, and we also study the limit of $U(t,m)-\lambda t$.
 
\subsection{Long-time behavior of HJ equations in finite dimensions}
Assuming that $H$ is independent of its measure variable, that is $H(x,p,m)=H_0(x,p)$ for some $H_0\colon \T^d\times \R^d\to \R$, and setting $m=\delta_x,\; x\in\T^d,$ in \eqref{HJB0}, we formally get the equation
\be\label{HJfd}
\begin{cases}
\partial_tu+H_0(x,Du)=0,& (t,x)\in(0,+\infty)\times\T^d,\\
u(0,x)=G(\delta_x)=:g(x),&x\in \T^d.
\end{cases}
\ee
The long-time behavior of $u$ solving \eqref{HJfd} in the (classical \cite{crandall1992user}) viscosity sense is tied to the ergodic problem
\be\label{EHJfd}
\lambda_0+H_0(x,Dv)=0,\;\;\text{in }\T^d.
\ee
This ergodic problem was introduced by Lions-Papanicolaou-Varadhan \cite{lions1987homogenization} in their study of periodic homogenization. They showed that, if $H_0$ is continuous and coercive in $p$\footnote{That is, $H_0(x,p)\to+\infty$ as $|p|\to+\infty$, uniformly in $x\in\T^d$. In particular, Hamiltonians with quadratic growth in $p$, as in this paper, are coercive.}, there exists a unique constant $\lambda_0$, called the ergodic constant, for which \eqref{EHJfd} has a viscosity solution $v\colon\T^d\to \R$. Since $v(x)+\lambda_0t$ solves the first equation in \eqref{HJfd}, a comparison argument shows that $u(t,x)-\lambda_0t$ is bounded, so that $u(t,x)/t\xrightarrow{t\to\infty}\lambda_0$ uniformly. The convergence of $u(t,x)-\lambda_0t$ itself is a much more delicate question, which requires additional structure on $H_0$.
\vspace{1mm}
 
\nd
For convex Hamiltonians, the study of this convergence was initiated by Namah-Roquejoffre \cite{namah1999remarks}, under a particular structure on $H_0$, and by Fathi \cite{fathi1998convergence}, who treated strictly convex and superlinear Hamiltonians via weak KAM theory, that is, through the dynamics of the associated Lagrangian flow; see also \cite{roquejoffre2001convergence}. Davini-Siconolfi \cite{davini2006generalized} identified the limit of $u(t,x)-\lambda_0 t$ through a Lax-type formula for strictly convex and coercive Hamiltonians, and obtained convergence results for Hamiltonians that are convex but not necessarily strictly convex. 
\vspace{1mm}
 
\nd
Most relevant to the present work is the approach of Barles-Souganidis \cite{barles2000large}, which does not rely on optimal control representation formulas and applies to non-convex Hamiltonians. They assume that $H_0$ is continuous and coercive, that the constant $\lambda_0$ in \eqref{EHJfd} is $0$, and a strict convexity-type assumption (SCA for short) holds:\\[0.3ex]
\null \hspace{0.5cm}For any $\eta>0$ small enough, there exists $\psi_\eta>0$ such that, for all $\mu\in (0,1]$ and $x\in\T^d$, $p,q\in\R^d$,
\be\label{SCAfd}\tag{SCA}
\left.\begin{aligned}
&H_0(x,p+q)\ge\eta\\
&H_0(x,q)\le0
\end{aligned}\right\}
\implies
\mu H_0\big(x,\frac{p}{\mu}+q\big)\ge H_0(x,p+q)+\psi_\eta(1-\mu).
\ee
Under these assumptions, it was proved that $u(t,\cdot)-\lambda_0t=u(t,\cdot)$ converges uniformly, as $t\to+\infty$, to a solution of \eqref{EHJfd}. We note that the key consequence of \eqref{SCAfd} that allows for the uniform convergence is an asymptotic monotonicity property for $u$, which formally implies that the negative part of $\partial_t(u-\lambda_0t)$ vanishes as $t\to+\infty$.

\vspace{1mm}

\nd
It is worth mentioning that all of the aforementioned results are proved on a compact state space (the torus or, more generally, a compact manifold). On the whole space $\R^d$, the convergence of $u(t,x)-\lambda_0t$ may fail even for convex Hamiltonians \cite{barles2000some,barles2006ergodic}; we refer to \cite{barles2006ergodic,ishii2008asymptotic} for positive results under additional assumptions. 
We finally refer to \cite{barles2012pde} for boundary-value problems with non-convex Hamiltonians and to \cite{tran2021hamilton} for a more extensive discussion on the literature.

\subsection{The problem in infinite dimensions and our results}
The literature on long-time behavior of HJ equations in the Wasserstein space is mostly based on its connection with the long-time and ergodic behavior of stochastic optimal control (or mean field control) problems or mean field games. Related works on the topic include \cite{bayraktar2025ergodicity,bayraktar2024infinite,cardaliaguet2020weak} for mean field control, \cite{cardaliaguet2012long,cardaliaguet2019long,cardaliaguet2013long} for mean field games, \cite{bayraktar2023solvability} for both, and the references therein. For linear-quadratic mean field control with common noise, \cite{bayraktarjian2026common} studies uniform-in-time convergence and turnpike properties in relation to an ergodic control problem.
\vspace{1mm}

\nd
Directly related to our setup are the works of Gangbo-Tudorascu \cite{gangbo2010lagrangian,gangbo2014weak}, who developed a weak KAM theory on $\mathcal{P}_2(\T^d)$ for mechanical Hamiltonians\footnote{Mechanical Hamiltonians have the form $H(x,p,m)=\frac12|p+c|^2+K(m)$, for some $K\colon \mathcal{P}_2\to \R$ and $c\in \R^d$.}, and the more recent work of Fuhrman-Rud\`a \cite{fuhrman2025ergodic}, who considered  mean field control problems (recall \eqref{MFCvalue}) on $\mathcal{P}_2(\R^d)$ with dissipative dynamics\footnote{We note that in the work \cite{fuhrman2025ergodic} the dynamics can also be influenced by an idiosyncratic noise; that is the state dynamics are allowed to take the form of a controlled diffusion process with a possibly non-trivial noise term:
$$dX_t=b(X_t,\alpha_t)dt+\sigma dB_t,\text{ where the dissipativity condition is }(b(x,\alpha)-b(y,\alpha))\cdot (x-y)\le -c|x-y|^2,\; x,y,\alpha\in \R^d,$$
for some constant $c>0$, for the drift $b\colon \R^d\times\R^d\to \R^d$ and the diffusion matrix $\sigma$. In our case, $b(x,\alpha)=\alpha$, hence this dissipativity condition does not hold.}. Both approaches rely on the control-theoretic or variational structure of the problem, and therefore do not apply when the Hamiltonian $H$ is not convex. However, we note that in order to study the long-time behavior, it was established that the infinite-dimensional analogue of \eqref{EHJfd} is the ergodic HJ equation
\be\label{EHJB}
\lambda+\int_{\T^d}H\big(x,D_mV(m,x),m\big)\,m(dx)=0,\qquad m\in\mathcal{P}_2,
\ee
where both the constant $\lambda$ (ergodic constant) and the function $V\colon\mathcal{P}_2\to\R$ are unknown. Finally, to compare with our own result (see also Remark \ref{genrem}), in the non-compact setting of $\mathcal{P}_2(\R^d)$, \cite{fuhrman2025ergodic} showed that the value function $W^T$ of the finite-horizon mean field control problem satisfies $W^T(0,m)-\lambda T=O(1)$ and $|W^T(0,m)-\lambda T-V(m)|\le C(1+{\bf d}_2^2(m,\delta_0))$, for all $m\in \mathcal{P}_2(\R^d)$.

\vspace{1mm}

\nd
In this paper, we use a purely PDE-based approach to study the ergodic problem \eqref{EHJB} and the long time behavior of the $L$-viscosity solution of \eqref{HJB0}. We first show that, if $H$ is locally Lipschitz continuous, with a Lipschitz constant that grows at most linearly in $p$, and has quadratic growth in $p$, then there exists a unique constant $\lambda\in\R$ for which \eqref{EHJB} admits a ${\bf d}_2$-Lipschitz $L$-viscosity solution (we refer to Definition \ref{evisc.}). Our main result (Theorem~\ref{main1}) states that if, in addition, $H$ satisfies an integral version of \eqref{SCAfd}, which we call (i-SCA) (see section~\ref{assumptionsandresults}), but is not necessarily convex, and $G$ is ${\bf d}_2$-Lipschitz, then there exists a ${\bf d}_2$-Lipschitz $L$-viscosity solution $V$ of \eqref{EHJB} such that
$$\lim_{t\to+\infty}\sup_{m\in\mathcal{P}_2}\big|U(t,m)-\lambda t-V(m)\big|=0.$$
We refer to section \ref{assumptionsandresults} for precise statements of our assumptions and results. To the best of our knowledge, this is the first convergence result of this type for HJ equations in the Wasserstein space with non-convex Hamiltonians. We finally present an example (section \ref{counter}) where long-time behavior may fail if \eqref{HJB0} is posed in $\mathcal{P}_2(\R^d)$ even if $H$ is quadratic.

\subsection{Organization of the paper}
We complete the introduction with the notation. In section~\ref{assumptionsandresults} we state our assumptions and the main result, and we outline its proof. In section~\ref{viscositysol} we recall the definition of $L$-viscosity solutions, collect the well-posedness and regularity results for \eqref{HJB0} that we use, and introduce the ergodic problem. In section \ref{ex}, we present applications of our main results to mean field control problems and zero-sum differential games. Section~\ref{proofs} is devoted to the construction of the ergodic constant and to the proof of our main result (Theorem~\ref{main1}). In section~\ref{counter} we provide an example where long-time convergence of $U(t,m)-\lambda t$ fails if the compact space $\mathcal{P}_2(\T^d)$ is substituted with the non-compact $\mathcal{P}_2(\R^d)$. Finally, Appendix~\ref{appendix} contains the proofs of some technical results. We provide examples of Hamiltonians satisfying our assumptions and assumption (i-SCA) in Appendix \ref{ex1}.

\subsection{Notation}
$\T^d=\R^d/\mathbb{Z}^d$ is the $d$-dimensional torus. For $\bm x\in (\T^d)^N$, we write $\bm x=(x^1,\ldots,x^N)$ for its $d$-dimensional coordinates. We use the notation $\bm 0=(0,\ldots,0)\in (\T^d)^N$. For a function $F\colon (\T^d)^N\to \R$, $D_{x^i}F(\bm x)$ denotes its gradient with respect to the $i\in\{1,\ldots,N\}$ coordinate. We also denote by $\text{Id}$ the $d\times d$ identity matrix.
\vspace{1mm}

\nd
For a function $\phi\colon\R\to \R$, we denote by $\phi'$ its derivative, if it exists. If the variable is specified, we write $\partial_t\phi (t)$ for the derivative.
\vspace{1mm}

\nd
$\mathcal{P}_2:=\mathcal{P}_2(\T^d)$ is the space of probability measures over $\T^d$. For a Polish space $(X,d)$, we endow the space $\mathcal{P}_2(X)$ of probability measures with finite second moment with the $2$-Wasserstein distance ${\bf d}_2$. For $p\in \{1,2\}$ we define the $p$-Wasserstein distance as
$${\bf d}_p(\mu,\nu):=\inf_{\pi\in \Pi(\mu,\nu)}\left(\iint d(x,y)^p\pi(dx,dy)\right)^{1/p},$$
where $\Pi(\mu,\nu)$ is the space of probability measures over $X\times X$ with marginals $\mu,\nu\in \mathcal{P}_2$. Note that ${\bf d}_1\le {\bf d}_2$. The space $\mathcal{P}_2$ endowed with ${\bf d}_2$ is referred to as Wasserstein space and it is known to be compact (see \cite{villani2009optimal}). For $\bm x\in (\T^d)^N$, we denote by $m^N_{\bm x}:=\frac{1}{N}\sum_{i=1}^N\delta_{x^i}\in \mathcal{P}_2$ the empirical measure.
\vspace{1mm}

\nd
Throughout the paper $(\Omega,\mathcal{F},\mathbb{P})$ is a fixed non-atomic probability space. We consider $\mathcal{H}=L^2(\Omega;\R^d)$, which is a separable Hilbert space endowed with the inner product $\langle X,Y \rangle=\mathbb{E}[X\cdot Y]$, $X,Y\in\mathcal{H}$. Its norm is denoted by $\|X\|_2=\left(\mathbb{E}[|X|^2]\right)^{1/2}$. For $X\in \mathcal{H}$, $\mathcal{L}(X)$ denotes its law. For a function $\Phi\colon [0,+\infty)\times \mathcal{H}\to \R$, we write $\Phi\in C^{1,1}$ if both $\partial_t\Phi$ and the Fr\'echet derivative $D\Phi$ exist and are continuous.\\
Let $\pi\colon \R^d\to \T^d$ be the natural projection. For $\mu\in\mathcal{P}_2(\R^d)$, $\pi_\#\mu(A):=\mu(\pi^{-1}(A)),$ for any measurable set $A\subset \T^d$.
\vspace{1mm}

\nd
For the functions $H\colon\T^d\times\R^d\times\mathcal{P}_2\to \R$ and $G\colon \mathcal{P}_2\to \R$ we use the notation
\be\label{hats}
\hat{H}(X,P)=\mathbb{E}[H(\pi(X),P,\pi_\#\mathcal{L}(X))],\quad \hat{G}(X)=G(\pi_\#\mathcal{L}(X)),\;\; \text{for }X,P\in \mathcal{H}.
\ee

\section{Assumptions and Main results}\label{assumptionsandresults}

\subsection{Assumptions}
We make the following assumptions on $H,G$.
\vspace{1mm}

\nd
\textbf{(A1).} $H\colon \T^d\times \R^d\times\mathcal{P}_2\to \R$ satisfies the following for some constants $C_H,c_H,K_0,k_0>0$.
\begin{enumerate}
    \item For any $(x,p,m),(x',p',m')\in \T^d\times\R^d\times\mathcal{P}_2$
    $$|H(x,p,m)-H(x',p',m')|\le C_H(1+|p|+|p'|)(|x-x'|+|p-p'|+{\bf d}_2(m,m')).$$
    \item  For any $(x,p,m)\in \T^d\times\R^d\times\mathcal{P}_2$
    $$c_H|p|^2-k_0\le H(x,p,m)\le C_H|p|^2+K_0.$$
\end{enumerate}
\vspace{2mm}

\nd
With the above assumptions (A1)(1) and (A1)(2), we are able to establish that there exists a unique constant $\lambda$ (ergodic constant) such that \eqref{EHJB} has a ${\bf d}_2$-Lipschitz
$L$-viscosity solution (see Proposition \ref{ergodicconstant}). We list the last assumption on $H$, which now depends on $\lambda$:
\begin{enumerate}   
    % \item[(i-SCA)] For every $\eta>0$ there exists a constant $\psi_\eta>0$, depending only on $\eta,\nu,c_H,C_H$, such that the following holds for all $\mu\in(0,1]$ and $\theta\in[0,\eta/2]$:\\[0.5ex]
    % If $X,P,Q\in\mathcal{H}$ satisfy $\hat{H}(X,P+Q)+\lambda\ge \eta$ and $\hat{H}(X,Q)+\lambda\le \theta$, then
    %     \be \label{SCA}\tag{\text{i-SCA}}
    %     \mu \left(\hat{H}(X,P/\mu+Q)+\lambda\right)\ge \hat{H}(X,P+Q)+\lambda+(\psi_\eta-\theta)(1-\mu),
    %     \ee
    %     where we used notation \eqref{hats}.
        \item[(i-SCA)] There exists a constant $c_*\ge1$ such that, for every $\eta>0$, there exists a constant $\psi_\eta>0$ with the following property, for all $\mu\in(0,1]$ and $\theta\in[0,\eta/2]$:\\[0.5ex]
    if $X,P,Q\in\mathcal{H}$ satisfy $\hat{H}(X,P+Q)+\lambda\ge \eta$ and $\hat{H}(X,Q)+\lambda\le \theta$, then
        \be \label{SCA}\tag{\text{i-SCA}}
        \mu \left(\hat{H}(X,P/\mu+Q)+\lambda\right)\ge \hat{H}(X,P+Q)+\lambda+(\psi_\eta-c_*\theta)(1-\mu),
        \ee
        where we used notation \eqref{hats}.
\end{enumerate}
\vspace{2mm}

\nd
\textbf{(A2).} $G\colon \mathcal{P}_2\to \R$ is ${\bf d}_2$-Lipschitz with constant $C_G$.
\vspace{2mm}

\nd
We denote by
\begin{align*}
&\text{(A1)}:\;\;\text{the combination of Assumptions (A1)(1) and (A1)(2)},\\
&\text{ (A)}:\;\; \text{the combination of Assumptions (A1) and (A2).}
\end{align*}
%Assumptions (A1) is unaffected if $H$ is translated by $\lambda$, hence once for the rest of the paper we assume that $\lambda=0$.

\begin{remark}[On the Assumptions]\label{genrem}
    (i) Assumption (A) is enough for us to establish a unique $L$-viscosity solution $U\colon [0,+\infty)\times\mathcal{P}_2\to \R$ of \eqref{HJB0}
    such that $U(t,\cdot)$ is ${\bf d}_2$-Lipschitz uniformly in $t\ge0$. The latter is obtained due to the quadratic growth in (A1)(2) (see Propositions \ref{VNLip}, \ref{ULip}).
    Under (A) we also obtain a unique ergodic constant $\lambda$ in \eqref{EHJB}.\\[0.3ex]
    (ii) Assumption (i-SCA) allows us to show that $U(t,m)-\lambda t$ is asymptotically monotone in $t$ (see Proposition \ref{AMP}) and gives us its uniform convergence.
    Comparing our setup with \cite{fuhrman2025ergodic}, our problem is posed in a compact space $\mathcal{P}_2(\T^d)$ and we can accommodate non-convex Hamiltonians $H$ satisfying (i-SCA) (see the next remark), while \cite{fuhrman2025ergodic}
    studies the mean field control problem \eqref{MFCvalue} in $\mathcal{P}_2(\R^d)$ with a dissipativity condition on the drift (to account for the lack of compactness) and obtains a bound of $U(t,m)-\lambda t-V(m),\; m\in\mathcal{P}_2(\R^d)$, for some $V\colon \mathcal{P}_2\to \R$ satisfying \eqref{EHJB} in the appropriate sense.
\end{remark}

\begin{remark}\label{asrem}
Two classes of examples of functions $H$ satisfying Assumption (i-SCA) are:\\[0.7ex]
    (i) Strongly convex Hamiltonians; that is, there exists $\nu>0$ such that for all $(x,m)\in\T^d\times\mathcal{P}_2$ the map $p\mapsto H(x,p,m)-\frac{\nu}{2}|p|^2$ is convex.\\[0.3ex]
    We refer to Lemma \ref{lem:SCA} for a proof of this fact. In Lemma \ref{lem:SCA} we also show that elements of this class of Hamiltonians satisfy \eqref{SCAfd} from \cite{barles2000large}, if they are independent of the measure variable.\\[0.5ex]
    (ii) Our second class of examples consists of Hamiltonians that are not necessarily convex. Let $\rho,s_0>0$ with $\rho<s_0/6$, and $f\in C^1(\T^d)$ so that $\min f<0<1<\max f$. We consider
    \be\label{costs}
    l(m):=\int_{\T^d}f(x) m(dx),\;\; F(m):= (l(m)-1)_+,\;\;\text{and }\;\;\chi(m)=\min\{ 1, (1-l(m))_+\}.
    \ee
    The Hamiltonian $H(x,p,m)=\frac{1}{2}|p|^2-\rho \chi(m)(|p|-s_0)_+-F(m)$ satisfies our Assumptions and is not convex in $p$. We prove this fact in Appendix \ref{ex1} (Lemma \ref{lem:gameISCA}). We also provide a connection of this Hamiltonian to a zero-sum differential game in section \ref{ex}.
\end{remark}

\subsection{Main result} The main result of this paper is summarized in the following theorem.

\begin{theorem}\label{main1}
    Assume that $H,G$ satisfy (A1)-(A2). There exists a unique $L$-viscosity solution $U\colon [0,+\infty)\times \mathcal{P}_2\to \R$ of \eqref{HJB0} and a unique constant $\lambda\in \R$ so that \eqref{EHJB} has a ${\bf d}_2$-Lipschitz $L$-viscosity solution.
    Furthermore, if Assumption (i-SCA) also holds, there exists a ${\bf d}_2$-Lipschitz $L$-viscosity solution $V\colon \mathcal{P}_2\to \R$ of \eqref{EHJB} such that
    \be\label{LTB}
    \lim_{t\to+\infty}\sup_{m\in\mathcal{P}_2}\left|U(t,m)-\lambda t- V(m)\right|=0.
    \ee
\end{theorem}

\nd
We provide some applications of the above results to Mean Field Control problems and to Zero-sum Differential Games in section \ref{ex}.

\subsection{Outline of the proof} Our starting point is that, under (A), \eqref{HJB0} satisfies a comparison principle between $L$-viscosity subsolutions and supersolutions and has a unique solution $U\colon [0,+\infty)\times \mathcal{P}_2\to \R$ such that $\mathcal{P}_2\ni m\mapsto U(t,m)$ is ${\bf d}_2$-Lipschitz with a constant independent of $t\ge 0$. This is achieved because of the quadratic growth of $H$ (Assumption (A1)(2)). In Proposition \ref{ergodicconstant}, we show that there exists a unique constant $\lambda\in \R$ such that $U(t,m)-\lambda t$ is a bounded function and such that \eqref{EHJB} (with constant $\lambda$) has an $L$-viscosity solution. 
In order to establish the long-time behavior \eqref{LTB}, we assume (i-SCA) on top of (A1)-(A2). We observe that, since Assumptions (A1) and (i-SCA) are invariant under translation by $\lambda$, we may assume without loss of generality that $\lambda=0$.\\[0.5ex]
Our argument for \eqref{LTB} is to show that the function $t\mapsto U(t,m)$ exhibits the asymptotic monotonicity property\footnote{At a formal level, this means that for all $m\in \mathcal{P}_2$, $\left( \partial_tU(t,m) \right)_-\to 0$, as $t$ becomes large.}; something first observed by Barles-Souganidis \cite{barles2000large} for \eqref{HJfd}.

\begin{proposition}[Asymptotic Monotonicity Property]\label{AMP}
        Under the Assumptions of Theorem \ref{main1}, let $U\colon [0,+\infty)\times\mathcal{P}_2\to \R$ be an $L$-viscosity solution of \eqref{HJB0} and assume without loss of generality that $\lambda=0$. For any $\eta >0$, there exists a bounded $\delta_\eta: [0,+\infty)\to [0,+\infty)$ such that 
        \begin{align}
            \delta_\eta(s)\xrightarrow{s\to +\infty}0\quad\text{and}\quad U(s,m)-U(t,m)+\eta (s-t)\le \delta_\eta(s),\text{ for all } m\in\mathcal{P}_2 \text{ and } t\ge s\ge0.\label{asmonot}
        \end{align}
        %\text{for all } $m\in\mathcal{P}_2$ \text{ and } $t\ge s\ge0$.
    \end{proposition}

\nd
We present the proof of the Asymptotic Monotonicity Property in section \ref{4.2}. Assuming that it holds, the proof of \eqref{LTB} proceeds in the following 3 steps.
\vspace{2mm}

\begin{enumerate}
    \item[Step 1.] Since $(U(t,\cdot))_{t\ge 0}$ is bounded and uniformly Lipschitz, there exists a sequence $(T_n)_{n\in \mathbb{N}}\subset [0,+\infty)$ with $T_n\xrightarrow{n\to +\infty} +\infty$, such that $(U(T_n,\cdot))_{n\in \mathbb{N}}$ converges uniformly in $\mathcal{P}_2$ as $n\to +\infty$. By the comparison principle, we may then establish
    % For $n,k\in\mathbb{N}$, by the comparison principle we may write
    % $$|U(T_n+t,m)-U(T_k+t,m)|\le \sup_{m\in\mathcal{P}_2}\{|U(T_n,m)-U(T_k,m)|\},\text{ for all }(t,m)\in (0,+\infty)\times\mathcal{P}_2.$$
    that $(U(T_n+\cdot,\cdot))_{n\in\mathbb{N}}$ is Cauchy in $C((0,+\infty)\times\mathcal{P}_2)$ and hence it converges to a function $U^\infty\in C((0,+\infty)\times\mathcal{P}_2)$. The function $U^\infty$ in fact satisfies the same Lipschitz estimate as $U$. 
    \vspace{1mm}

    \item[Step 2.] By Proposition \ref{AMP}, we have for any $m\in \mathcal{P}_2$, $t\ge s\ge 0$,  $n\in\mathbb{N}$ and $\eta>0$ 
    $$U(s+T_n,m)-U(t+T_n,m)+\eta(s-t)\le \delta_\eta(s+T_n).$$
    Sending $n\to +\infty$ and then $\eta\to 0$, this gives $U^\infty(s,m)\le U^\infty(t,m)$ for all $m\in\mathcal{P}_2$ and $t\ge s\ge 0$. Since $m\mapsto U^\infty(t,m)$ is bounded and equilipschitz, we conclude that $U^\infty(t,m)\xrightarrow{t\to+\infty}V(m)$ uniformly (increasing limit), where $V\colon \mathcal{P}_2\to \R$ is ${\bf d}_2$-Lipschitz.
    \vspace{1mm}
    
    \item[Step 3.] The uniform convergence of $U(t+T_n,m)$ to $U^\infty(t,m)$ implies $U(t,m)\xrightarrow{t\to +\infty}V(m)$ uniformly. We complete the proof by showing that $V$ is an $L$-viscosity solution of \eqref{EHJB} by a stability argument on the equation satisfied by $(U(\cdot+k,\cdot))_{k\in \mathbb{N}}$.
\end{enumerate}

\section{Viscosity solutions}\label{viscositysol}

\nd
In this section, consider for $T\in (0,+\infty)\cup \{+\infty\}$ the HJ equation
\be\label{HJB}
\begin{cases}
\partial_tU+\int_{\T^d}H(x,D_mU(t,m,x),m)m(dx)=0,&(t,m)\in (0,T)\times \mathcal{P}_2,\\
U(0,m)=G(m),& m\in\mathcal{P}_2.
\end{cases}
\ee

\subsection{Definitions}
Viscosity solutions of \eqref{HJB0} (or the ergodic problem) are understood in the lifted sense ($L$-viscosity solutions). We define the ``lift'' of a function defined in $\mathcal{P}_2$. Let $\pi\colon \R^d\to \T^d$ be the natural projection.

\begin{definition}
    For a given $\Phi:\mathcal{P}_2\rightarrow \R$, its ``lift'' $\hat{\Phi}:\mathcal{H}\rightarrow \R$ is defined by
    $\hat{\Phi}(X):=\Phi(\pi_\#\mathcal{L}(X)).$
\end{definition}

\begin{remark}
Recall that for a function $\hat{\Phi}:\mathcal{H}\rightarrow \R$ its gradient $D \hat{\Phi}(X)\in \mathcal{H}$ satisfies
$$\hat{\Phi}(X+Y)=\hat{\Phi}(X)+\langle D\hat{\Phi}(X), Y\rangle+o(\|Y\|_2),\;\; X,Y\in\mathcal{H}.$$
If $\Phi$ is smooth enough (see \cite{carmona2018probabilistic}), the gradient $D\hat{\Phi}$ is related to the Wasserstein derivative $D_m\Phi$ in the following way
$$D_m\Phi(\pi_\#\mathcal{L}(X),\pi(X))=D\hat{\Phi}(X),\;\; X\in\mathcal{H}.$$
\end{remark}
\vspace{1mm}

\nd
With this in mind, we notice that, at least formally, if $U$ satisfies \eqref{HJB}, then its ``lift'' $\hat{U}:[0,T)\times\mathcal{H}\rightarrow \R$ satisfies
\be\label{LHJB}\tag{$\mathrm{HJB}_\mathcal{H}$}
\begin{cases}
    \partial_t \hat{U} +\hat{H}(X,D\hat{U}(t,X))=0, &(t,X)\in (0,T)\times \mathcal{H},\\
    \hat{U}(0,X)=\hat{G}(X), &X\in\mathcal{H}.
\end{cases}
\ee

\begin{definition}\label{L-visc}
    Let $T\in (0,+\infty)\cup\{+\infty\}$. An upper semicontinuous (resp. lower semicontinuous) function $V:[0,T)\times \mathcal{H}\rightarrow \R$ is a viscosity subsolution (resp. supersolution) of \eqref{LHJB} if
    \begin{itemize}
        \item $V(0,X)\le (\text{resp. }\ge)\; \hat{G}(X)$;
        \vspace{1mm}
        \item For any $\Phi\in C^{1,1}([0,T)\times\mathcal{H})$ and each $(t_0,X_0)\in (0,T)\times \mathcal{H}$ such that $V(t_0,X_0)-\Phi(t_0,X_0)$ has a local maximum (resp. local minimum), we have
        $$\partial_t \Phi(t_0,X_0)+\hat{H}(X_0,D\Phi(t_0,X_0))\le (\text{resp. }\ge)\; 0.$$
    \end{itemize}
    $V$ is a viscosity solution of \eqref{LHJB} if it is a viscosity subsolution and a viscosity supersolution.
\end{definition}

\begin{definition}\label{visc.}
    Let $T\in (0,+\infty)\cup\{+\infty\}$. A function $U\colon [0,T)\times\mathcal{P}_2\rightarrow \R$ is a viscosity solution/subsolution/supersolution of \eqref{HJB} if its lift $\hat{U}$ is a viscosity solution/subsolution/supersolution of \eqref{LHJB}.
\end{definition}

\subsection{Time-dependent problem}
As is established in the mean field control literature, an $L$-viscosity solution of \eqref{HJB0} arises as a limit (in the appropriate sense) of the solutions of the finite dimensional problems
\be\label{NHJB}
\begin{cases}
    \partial_t V^N +\frac{1}{N}\sum_{i=1}^NH(x^i,ND_{x^i}V^N,m^N_{\bm x})=0,&(t,\bm x)\in[0,+\infty)\times(\T^d)^N,\\
    V^N(0,\bm x)= G(m^N_{\bm x}),& \bm x\in (\T^d)^N,
\end{cases}
\ee
where $N\in \mathbb{N}$. We have the following two Propositions.
\vspace{2mm}

\begin{proposition}\label{VNLip}
    Assume that (A) holds and let $V^N$ be the unique viscosity solution of \eqref{NHJB}. Then there exists a constant $C_0>0$ depending only on $H,G$ such that for every $N\in \mathbb{N}$ we have
    $$|V^N(t,\bm x)-V^N(s,\bm y)|\le C_0\left( {\bf d}_2(m^N_{\bm x},m^N_{\bm y})+|t-s|\right),\;\text{for all }t,s\ge0 \text{ and } \bm x,\bm y\in (\T^d)^N.$$
\end{proposition}
\vspace{1mm}

\begin{proof}
    We give the proof in Appendix \ref{appendix}.
\end{proof}

\nd
As a consequence, we obtain:
\begin{proposition}\label{ULip}
    Assume (A). There exists a constant $C_0>0$ depending only on $H,G$ and $U\colon [0,+\infty)\times\mathcal{P}_2\to \R$ an $L$-viscosity solution of \eqref{HJB0} such that 
    $$|U(t,m)-U(s,n)|\le C_0({\bf d}_2(m,n)+|t-s|),\;\;\text{for all }t,s\ge 0,\;m,n\in\mathcal{P}_2.$$
\end{proposition}

\begin{proof}
    We give the proof in Appendix \ref{appendix}.
\end{proof}

\nd
The following theorem states that \eqref{HJB0} admits a comparison principle for continuous $L$-viscosity subsolutions and continuous supersolutions. Therefore, $U$ from Proposition \ref{ULip} is the unique $L$-viscosity solution of \eqref{HJB0}.

\begin{theorem}\label{comp}
    Assume that $H$ satisfies Assumption (A1) and consider $U_1,U_2\colon [0,+\infty)\times\mathcal{P}_2\to \R$ be a continuous $L$-viscosity subsolution and a continuous $L$-viscosity supersolution of the first equation in \eqref{HJB}, respectively, such that $U_1(0,\cdot)\le U_2(0,\cdot)$. Then, $U_1(t,m)\le U_2(t,m)$, for all $(t,m)\in [0,+\infty)\times\mathcal{P}_2$. In particular, there exists at most one $L$-viscosity solution of \eqref{HJB0}.
\end{theorem}

\begin{proof}
    We use \cite[Theorem 7.4]{gangbo2021finite} that provides us with a comparison principle between bounded subsolutions and bounded supersolutions of \eqref{HJB} if $T<\infty$. We note that for any $T<\infty$, the restrictions $U_1^T,U_2^T\colon [0,T]\times\mathcal{P}_2\to \R$ of $U_1,U_2$ are an $L$-viscosity subsolution (respectively supersolution) of \eqref{HJB}. Since $U_1^T,U_2^T$ are continuous on the compact set $[0,T]\times \mathcal{P}_2$, they are also bounded. Thus, $U_1^T\le U_2^T$, which gives the desired $U_1\le U_2$.
\end{proof}

\begin{remark}\label{MFCrem}
    If on top of (A1), $H(x,p,m)$ is assumed to be convex in the $p$ variable, then it follows by applying the dynamic programming principle that the value function $W^T(t,m)$ from \eqref{MFCvalue} satisfies the HJ equation 
    $$\begin{cases}
    -\partial_tU+\int_{\T^d}H\big(x,D_mU(t,m,x),m\big)\,m(dx)=0,& (t,m)\in (0,T)\times\mathcal{P}_2,\\
    U(T,m)=G(m), &m\in \mathcal{P}_2,
\end{cases}$$
in the $L$-viscosity sense. Thus, by the proof of Theorem \ref{comp}, $U(T-t,m)=W^T(t,m)$, for all $(t,m)\in [0,T]\times \mathcal{P}_2$, where $U$ is the function from Proposition \ref{ULip}.
\end{remark}

% \begin{proof}
%     Fix $T>0$. We reverse the time and we define $U_1^T,U_2^T\colon [0,T]\times\mathcal{P}_2\to \R$ with $U^T(t,m):=U(T-t,m)$, for all $(t,m)\in [0,T]\times \mathcal{P}_2$. Then, $U^T_1,U_2^T$ are $L$-viscosity subsolution and $L$-viscosity supersolution of the backwards HJB with $G$ as terminal function at time $T$
%     $$\begin{cases}
%         -\partial_tW+\int H(x,D_mW(t,m,x),m)m(dx)=0,&(t,m)\in [0,T)\times\mathcal{P}_2,\\
%         W(T,m)=G(m),& m\in \mathcal{P}_2.
%     \end{cases}$$
%     The comparison principle for these problems \cite[Theorem 7.4]{gangbo2021finite} yields $U_1^T\le U_2^T$. Since $T>0$ was chosen arbitrarily, we conclude that $U_1\le U_2$ and the uniqueness of $L$-viscosity solutions follows.
% \end{proof}

\subsection{Ergodic problem}
In this subsection we provide the definition of $L$-viscosity solutions of \eqref{EHJB}.

\begin{definition}\label{evisc.}
    (i) Let $\lambda\in \R$. An upper (resp. lower) semicontinuous function $V\colon \mathcal{P}_2\to \R$ is an $L$-viscosity sub(resp. super)solution of \eqref{EHJB} if its lift $\hat{V}\colon \mathcal{H}\to \R$ is a viscosity sub(resp. super)solution of 
    $$\lambda+\hat{H}(X,D\hat{V}(X))=0,\;\;X\in \mathcal{H},$$
    where $\hat{H}$ was defined in \eqref{hats}.\\[0.5ex]
    (ii) Let $\lambda\in \R$. A continuous function $V\colon \mathcal{P}_2\to \R$ is an $L$-viscosity solution of \eqref{EHJB} if it is both an $L$-viscosity subsolution and supersolution of \eqref{EHJB}. 
    %In this case, we also say that $(V,\lambda)$ is a solution of \eqref{EHJB}.
\end{definition}

%\nd
%We focus on all constants $\lambda$ so that \eqref{EHJB} has a ${\bf d}_2$-Lipschitz $L$-viscosity solution. If such a constant exists, it is called  \textit{ergodic constant}. In section \ref{proofs}, we show that there exists a unique ergodic constant.

\section{Applications}\label{ex}

\nd
In this section we present some applications of our main result to Mean Field Control problems and to Zero-Sum Differential Games.

\subsection{Mean Field Control problems} Assume that there exists $\nu>0$ such that the map $\R^d\ni p\mapsto H(x,p,m)-\frac{\nu}{2}|p|^2$ is convex for any $(x,m)\in \T^d\times \mathcal{P}_2$. By Remark \ref{MFCrem}, for any $T>0$, we know that the value of the Mean Field Control problem \eqref{MFCvalue} satisfies $W^T(t,m)=U(T-t,m)$, for any $(t,m)\in [0,T]\times\mathcal{P}_2$, where $U$ is the unique $L$-viscosity solution of \eqref{HJB0}. Then, Theorem \ref{main1} in particular implies that
$$W^T(0,m)-\lambda T-V(m)\xrightarrow{T\to +\infty}0,\;\;\text{uniformly in }m\in \mathcal{P}_2,$$
for the ergodic constant $\lambda$ and some function $V\colon \mathcal{P}_2\to \R$ satisfying \eqref{EHJB} in the $L$-viscosity sense.

\subsection{Zero-sum differential games} Fix $T>0$, $(t,m)\in[0,T]\times\mathcal{P}_2$ and recall the functions $l, F,\chi$ from \eqref{costs}. We consider the running cost $\ell_0\colon \R^d\times \R^d\times \mathcal{P}_2\to \R$ with
$$\ell_0(\alpha,\beta ,m):=\frac12|\alpha|^2-s_0\,\chi(m)\,|\beta|+F(m).$$
For this section only, we make the following assumption
\vspace{2mm}

\nd
\textbf{Assumption}:\;\; (A3)\hspace{1cm} $G\colon \mathcal{P}_2\to \R$ is ${\bf d}_1$-Lipschitz with constant $L_G$.
\vspace{2mm}

\nd
We note that Assumption (A3) implies Assumption (A2).
\vspace{2mm}

\nd
We consider a filtered probability space $(\widetilde{\Omega},\widetilde{\mathcal{F}},\widetilde{\mathbb{F}},\widetilde{\mathbb{P}})$ which is rich enough and accommodates a $d$-dimensional Brownian motion $B=(B_t)_{0\le t\le T}$. We also assume that $\widetilde{\mathcal{F}}_0$ is atomless, and $\widetilde{\mathbb{F}}$ is the augmentation of the filtration generated by $\widetilde{\mathcal{F}}_0$ and $B$. 
\vspace{1mm}

\nd
The minimizing player chooses an $\widetilde{\mathbb{F}}$-progressively measurable $\alpha\colon [t,T]\times \widetilde{\Omega}\to \R^d$ with $\int_t^T\widetilde{\mathbb{E}}[|\alpha_s|^2]ds<\infty$, while the maximizing player chooses an $\widetilde{\mathbb{F}}$-progressively measurable $\beta\colon [t,T]\times \widetilde{\Omega}\to \R^d$ with $|\beta_s|\le\rho$ a.s. for all $s\in [t,T]$.
% The minimizing player chooses a
% measurable velocity field $\alpha:[t,T]\times\mathbb{T}^d\to\mathbb{R}^d$, while the maximizing player chooses a measurable field $\beta:[t,T]\times\mathbb{T}^d\to\mathbb{R}^d$ with $|b_s(x)|\le\rho$ for all $(s,x)$. 
We denote the sets of admissible controls by $\mathcal{A}_t$ and
$\mathcal{B}_t$. For $\alpha\in \mathcal{A}_t$ and $\beta\in \mathcal{B}_t$, the
population evolves according to
\begin{equation}\label{game:dyn}
dX_s=(\alpha_s+\chi(m_s)\,\beta_s)ds,\qquad
m_s=\mathcal{L}(X_s),\qquad \mathcal{L}(X_t)=m,
\end{equation}
and the payoff is
\begin{equation}\label{game:payoff}
J(t,m;\alpha,\beta):=\widetilde{\mathbb{E}}\left[\int_t^T
\ell_0\big(\alpha_s,\beta_s,m_s\big)\,ds+G(m_T)\right].
\end{equation}
From the particular structure of our $\ell_0$, we see that the minimizer pays a kinetic cost and the running cost $F$. The maximizer gains from a large payoff, but pays a cost of $s_0$ per unit of his/her actual speed $\chi(m_s)|\beta_s|$.
\vspace{1mm}

\nd
We now define the sets of non-anticipative strategies (see \cite{elliott1972existence}).

\begin{definition}\label{nonanti}
    A map
$\mathfrak{a}:\mathcal{B}_t\to\mathcal{A}_t$ is a non-anticipative
strategy for the minimizer if, for every $s\in[t,T]$ and all
$\beta,\beta'\in\mathcal{B}_t$,
\[
\widetilde{\mathbb{P}}\left(\beta_\tau=\beta'_\tau\ \text{ for a.e. }\tau\in[t,s]\right)=1
\Longrightarrow \widetilde{\mathbb{P}}\left(
\mathfrak{a}[\beta]_\tau=\mathfrak{a}[\beta']_\tau\ \text{ for a.e. }\tau\in[t,s]\right)=1.
\]
Non-anticipative strategies $\mathfrak{b}:\mathcal{A}_t\to\mathcal{B}_t$ for
the maximizer are defined in the same way. We denote the two sets of
strategies by $\mathfrak{S}^{\min}_t$ and $\mathfrak{S}^{\max}_t$.
\end{definition}

\nd
The lower and upper values of this zero-sum game are
\begin{align}\label{valuef}
W^-(t,m):=\inf_{\mathfrak{a}\in\mathfrak{S}^{\min}_t}\ \sup_{\beta\in\mathcal{B}_t}
J\big(t,m;\mathfrak{a}[\beta],\beta\big)\quad\text{and}\quad
W^+(t,m):=\sup_{\mathfrak{b}\in\mathfrak{S}^{\max}_t}\ \inf_{\alpha\in\mathcal{A}_t}
J\big(t,m;\alpha,\mathfrak{b}[\alpha]\big);
\end{align}
that is, each player can choose a strategy to react to the opponent's current action. In $W^-$ the
advantage belongs to the minimizer, and in $W^+$ to the maximizer. 
%We note that $W^-\le W^+$.

\subsubsection{Dynamic programming principle and associated HJ equations} We first show that in the above zero-sum game, the control set $\mathcal{A}_t$ can be restricted to controls $\alpha\in \mathcal{A}_t$ such that $|\alpha|\le R$, for some $R>0$ large enough. This will allow us to use the dynamic programming principle from \cite[Corollary 4.1]{cosso2019zero} and its consequent viscosity properties of the value functions \eqref{valuef}. In what follows, for $R>0$, we denote by $\pi_R\colon \R^d\to \R^d$ the projection
$$\pi_R(p)=\begin{cases}
    p,&\text{if }|p|\le R,\\
    R\frac{p}{|p|},&\text{if }|p|>R,
\end{cases}$$
and by $\mathcal{A}_t^R$ the set of admissible controls $\alpha\in\mathcal{A}_t$ with $|\alpha|\le R$.

\begin{lemma}\label{lem:projection-game}
Assume that $G$ satisfies Assumption (A3). Let $L_f$ be the Lipschitz constant of $F,\chi,l$ from Remark \ref{asrem}(ii) and set $C_T:=e^{\rho L_fT}\big(T\,L_f(1+s_0\rho)+L_G\big).$ Then, for every $R\ge C_T$, $t\in[0,T]$, $m\in\mathcal{P}_2$, $X\in L^2(\widetilde{\Omega},\widetilde{\mathbb{P}})$ which is $\widetilde{\mathcal{F}}_t$-measurable with $m=\mathcal{L}(X)$,
$\alpha\in\mathcal{A}_t$ and $\beta\in\mathcal{B}_t$,
$$J(t,m;\pi_R(\alpha),\beta)\le J(t,m;\alpha,\beta).$$
\end{lemma}

\begin{proof}
For $s\in [t,T]$, let $X_s$ and $X^R_s$ be the states in \eqref{game:dyn} starting from $X_t=X_t^R=X$ driven by $(\alpha,\beta)$ and
$(\pi_R(\alpha),\beta)$, with laws $m_s$ and $m^R_s$, and set
$$0\le D:=\widetilde{\mathbb{E}}\left[\int_t^T(|\alpha_s|-R)_+\,ds\right]<\infty.$$
Let $a\in \R^d$. We note that if $|a|>R$, then $\frac12(|a|^2-|\pi_R(a)|^2)=\frac12(|a|-R)(|a|+R)\ge R(|a|-R)$; if $|a|\le R$,
both sides vanish. Hence
\be\label{eq:1}
\widetilde{\mathbb{E}}\left[\int_t^T\frac12\big(|\pi_R(\alpha_s)|^2-|\alpha_s|^2\big)\,ds\right]\le-R\,D.
\ee
Since $|\alpha_r-\pi_R(\alpha_r)|=(|\alpha_r|-R)_+$
and $|\beta_r|\le\rho$,
$$|X_s-X^R_s|\le\int_t^s(|\alpha_r|-R)_+\,dr+\rho L_f\int_t^s{\bf d}_1(m_r,m^R_r)\,dr.$$
Let $e(r):=\widetilde{\mathbb{E}}[|X_r-X^R_r|]$, so that ${\bf d}_1(m_r,m^R_r)\le e(r)$. Then, taking expectations in the above gives
$e(s)\le D+\rho L_f\int_t^se(r)\,dr$, and Gr\"onwall's lemma implies
$${\bf d}_1(m_s,m^R_s)\le e(s)\le e^{\rho L_fT}D,\qquad s\in[t,T].$$
By the Lipschitz continuity of $F$, $\chi$ and $G$, and since $|\beta|\le\rho$,
\be\label{eq:2}
\widetilde{\mathbb{E}}\left[\int_t^T\Big(F(m^R_s)-F(m_s)+s_0|\beta_s|\big(\chi(m_s)-\chi(m^R_s)\big)\Big)\,ds\right]\le T\,L_f(1+s_0\rho)\,e^{\rho L_fT}D,
\ee
and
\be\label{eq:3}
G(m^R_T)-G(m_T)\le L_G{\bf d}_1(m^R_T,m_T)\le L_G\,e^{\rho L_fT}D.
\ee
Adding \eqref{eq:1}, \eqref{eq:2} and \eqref{eq:3} yields 
$$J(t,m;\pi_R(\alpha),\beta)-J(t,m;\alpha,\beta)\le(C_T-R)\,D\le0,$$
as desired.
\end{proof}

\nd
Let $W^-_R, W^+_R$ be the value functions from \eqref{valuef} if $\mathcal{A}_t$ is replaced by $\mathcal{A}_t^R$. For the non-anticipative strategies we use the notation $\mathfrak{S}_t^{\min,R},\;\mathfrak{S}_t^{\max,R}$ in this case. We show that if $R\ge C_T$, then $W^\pm_R=W^\pm$.

\begin{proposition}\label{prop:truncation-game}
Assume that $G$ satisfies Assumption (A3). Then, for every $R\ge C_T$,
$t\in[0,T]$ and $m\in\mathcal{P}_2$, 
$$W^-_R(t,m)=W^-(t,m)\quad\text{and}\quad W^+_R(t,m)=W^+(t,m).$$
\end{proposition}

\begin{proof}
We prove two inequalities for each value.
\vspace{1mm}

\nd
For $W^-_R=W^-$, we immediately notice that $W^-\le W^-_R$, since the infimum defining $W^-$ is over the larger set
$\mathfrak{S}^{\min}_t\supset\mathfrak{S}^{\min,R}_t$. For the opposite inequality, we observe that for $\mathfrak{a}\in\mathfrak{S}^{\min}_t$, the strategy
$\pi_R\circ\mathfrak{a}$ is non-anticipative and belongs to
$\mathfrak{S}^{\min,R}_t$. By Lemma~\ref{lem:projection-game},
$J(t,m;\pi_R(\mathfrak{a}[\beta]),\beta)\le J(t,m;\mathfrak{a}[\beta],\beta)$ for all $\beta\in \mathcal{B}_t$. Taking the supremum over $\beta$ and then the infimum over $\mathfrak{a}$, yields the opposite inequality.
\vspace{1mm}

\nd
We now show that $W^+_R=W^+$. Restricting $\mathfrak{b}\in\mathfrak{S}^{\max}_t$ to
$\mathcal{A}^R_t$ gives an element of $\mathfrak{S}^{\max,R}_t$, and an
infimum over $\mathcal{A}_t$ is at most the infimum over $\mathcal{A}^R_t$, thus $W^+\le W^+_R$. For the opposite inequality, let $\mathfrak{b}_R\in\mathfrak{S}^{\max,R}_t$, set $\mathfrak{b}[\alpha]:=\mathfrak{b}_R[\pi_R(\alpha)]$, which belongs to $\mathfrak{S}^{\max}_t$. By Lemma~\ref{lem:projection-game},
$$J(t,m;\alpha,\mathfrak{b}[\alpha])\ge J(t,m;\pi_R(\alpha),\mathfrak{b}_R[\pi_R(\alpha)])\ge\inf_{\alpha'\in\mathcal{A}^R_t}J(t,m;\alpha',\mathfrak{b}_R[\alpha']),$$ 
\vspace{-1mm}
for every $\alpha\in\mathcal{A}_t$. We take the infimum over $\alpha$ and then the
supremum over $\mathfrak{b}_R$ to discover $W^+\ge W^+_R$.
\end{proof}

\nd
We may now obtain the long time behavior of $W^\pm$ from \eqref{valuef}, assuming that (A3) is in place. We use the dynamic programming
principle \cite[Corollary 4.1]{cosso2019zero} to show that $W^\pm_R,\; R\ge C_T$ satisfies
\[\begin{cases}
-\partial_tW^\pm_R+\int_{\mathbb{T}^d}H^\pm_R\big(x,D_mW^\pm_R(t,m,x),m\big)\,m(dx)=0,& (t,m)\in [0,T)\times\mathcal{P}_2,\\
W^\pm_R(T,m)=G(m),& m\in \mathcal{P}_2,
\end{cases}
\]
with the matching sign, in the $L$-viscosity sense (see \cite[Section 5]{cosso2019zero} for a proof), where $L$-viscosity solutions for finite time HJ equations can be defined as in Definition \ref{visc.} and
\begin{align*}
H^-_R(x,p,m)&:=\inf_{|\beta|\le\rho}\ \sup_{\alpha\in\mathbb{R}^d,\;|\alpha|\le R}
\Big\{-\big(\alpha+\chi(m)\beta\big)\cdot p-\ell_0(\alpha,\beta,m)\Big\},\\
H^+_R(x,p,m)&:=\sup_{\alpha\in\mathbb{R}^d,\;|\alpha|\le R}\ \inf_{|\beta|\le\rho}
\Big\{-\big(\alpha+\chi(m)\beta\big)\cdot p-\ell_0(\alpha,\beta,m)\Big\}.
\end{align*}
It is easy to calculate that $H^\pm_R\xrightarrow{R\to +\infty}H$, where $H(x,p,m)=\frac{1}{2}|p|^2-\rho \chi(m)(|p|-s_0)_+-F(m)$ is the Hamiltonian from Remark \ref{asrem}(ii) and satisfies (A1) and (i-SCA) (see Lemma \ref{lem:gameISCA}). Since $W^\pm_R=W^\pm$ (with the matching sign) by Proposition \ref{prop:truncation-game}, we conclude that $W^+, W^-$ satisfy the HJ equation 
\[\begin{cases}
-\partial_tW^\pm+\int_{\mathbb{T}^d}H\big(x,D_mW^\pm(t,m,x),m\big)\,m(dx)=0,& (t,m)\in [0,T)\times\mathcal{P}_2,\\
W^\pm(T,m)=G(m),& m\in \mathcal{P}_2,
\end{cases}
\]
in the $L$-viscosity sense.
%Using our particular cost structure, we compute $H^+=H^-=H$ (Isaacs condition), where $H(x,p,m)=\frac{1}{2}|p|^2-\rho \chi(m)(|p|-s_0)_+-F(m)$. 
Therefore, due to the comparison principle, $W^-(t,m)=W^+(t,m)=U(T-t,m)$, where $U$ is the $L$-viscosity solution of \eqref{HJB0} with the same Hamiltonian $H$ provided by Proposition \ref{ULip}. Theorem \ref{main1} gives the long time behavior of $W^+(t,m)=W^-(t,m)$ as $T\to +\infty$.

\begin{remark}
    We expect that Theorem \ref{main1} can be applied to obtain the long-time behavior of the value functions of zero-sum games \eqref{valuef} with more general running costs $\ell_0$, provided that the Isaacs condition holds, the minimizer's controls can be truncated as in Lemma \ref{lem:projection-game}, and the (associated) Hamiltonian $H$ satisfies (A1) and (i-SCA). We do not pursue this direction in the present paper.
\end{remark}

\section{Proof of the main result}\label{proofs}

\nd
Before discussing the large time behavior of $U$ solving \eqref{HJB0}, we study its ergodic behavior. This leads to the construction of the ergodic constant $\lambda$.

\subsection{Construction of the ergodic constant and its uniqueness}

\begin{proposition}\label{ergodicconstant}
    Assume that $H$ satisfies (A1). Then there exists exactly one ergodic constant $\lambda\in \R$ such that \eqref{EHJB} has a ${\bf d}_2$-Lipschitz viscosity solution. Furthermore, if $U$ is the $L$-viscosity solution of \eqref{HJB0}, the map $[0,+\infty)\times \mathcal{P}_2\ni(t,m)\mapsto U(t,m)-\lambda t$ is bounded.
    %Furthermore, $\lambda$ is the unique constant such that for any $L$-viscosity solution $U$ of \eqref{HJB0} with $U(0,\cdot)$ ${\bf d}_1$-Lipschitz, $t\mapsto U(t,m)-\lambda t$ is uniformly bounded.
\end{proposition}

\begin{proof} 
We divide the proof in 2 steps.
\vspace{1mm}

\nd
\textit{Step 1.} (Existence of $\lambda$)\\
By classical results from finite dimensional HJ equations, for any $\delta>0$, there exists a unique $V^{N,\delta}=V^{N,\delta}(x^1,\ldots,x^N)\colon (\T^d)^N\to \R$, which is symmetric in its variables, viscosity solution of
$$\delta V^{N,\delta}(\bm x)+\frac{1}{N}\sum_{i=1}^NH(x^i,ND_{x^i}V^{N,\delta}(\bm x),m^N_{\bm x})=0,\;\;\bm x\in (\T^d)^N.$$
Furthermore, using the result of \cite{lions1987homogenization}, sending $\delta\to 0$ (along a subsequence), we obtain a constant $\lambda^N\in \R$ and a function $V^N\colon (\T^d)^N\to \R$ that satisfy
$$\delta V^{N,\delta}(\bm x)\xrightarrow{\delta\to 0}\lambda^N\quad\text{and}\quad V^{N,\delta}(\bm x)-V^{N,\delta}(\bm 0)\xrightarrow{\delta\to 0} V^N(\bm x),$$
uniformly in $\bm x\in (\T^d)^N$, so that 
\be\label{Nerg}\nonumber
    \lambda^N+\frac{1}{N}\sum_{i=1}^NH(x^i,ND_{x^i}V^N(\bm x),m^N_{\bm x})=0,\;\;\bm x\in (\T^d)^N,
\ee
is satisfied in the viscosity sense. By this construction, we note that $V^N(\bm 0)=0$ and $V^N(x^1,\ldots,x^N)$ is symmetric in its variables $x^1,\ldots,x^N$, because $V^{N,\delta}$ is symmetric in its variables.
\vspace{1mm}

\nd
Assuming, now, that $V^N$ achieves its maximum at $\bm x_0\in (\T^d)^N$, we discover, due to the quadratic growth of $H$ (Assumption (A1)(2)),
    $$0\ge \lambda^N+\frac{1}{N}\sum_{i=1}^NH(x_0^i,0,m^N_{\bm x_0})\ge \lambda^N -k_0.$$
    Hence the sequence $(\lambda^N)_{N\in\mathbb{N}}$ is bounded from above. A similar argument with the minimum of $V^N$ implies that $(\lambda^N)_{N\in\mathbb{N}}$ is bounded from below, hence $(\lambda^N)_{N\in\mathbb{N}}$ is bounded. Using the techniques from Proposition \ref{VNLip}
    %\cite[Proposition 2.5]{daudin2025error} 
    and the boundedness of $\lambda^N$, we also get
    %\eqref{Nerg} also implies that
    %$$N\sum_{i=1}^N\left|D_{x^i}V^N(\bm x)\right|^2\le C,\;\;\text{for a.e. }\bm x\in (\T^d)^N,$$
    %for some constant $C$ independent of $N$. This bound implies that $(V^N)_{N\in\mathbb{N}}$ satisfies
    $$|V^N(\bm x)-V^N(\bm y)|\le C{\bf d}_2(m^N_{\bm x},m^N_{\bm y}),$$ 
    for some constant $C$ independent of $N$.
    \vspace{1mm}

    \nd
    Now, since $V^N$ are symmetric, we define $U^N(m^N_{\bm x}):=V^N(\bm x)$. By the uniform ${\bf d}_2$-Lipschitz property of $U^N$, we can extend them to $\tilde{U}^N\colon \mathcal{P}_2\to \R$, which are also uniformly ${\bf d}_2$-Lipschitz. By Arzelà–Ascoli and Bolzano-Weierstrass we then obtain that there exists a ${\bf d}_2$-Lipschitz function $V\colon \mathcal{P}_2\to \R$ and a constant $\lambda$ such that
    $$\tilde{U}^N\xrightarrow{N\to +\infty}V\;\;\text{uniformly in}\;\mathcal{P}_2,\quad\text{and}\quad\lambda^N\xrightarrow{N\to +\infty}\lambda,$$
    up to a subsequence. The argument concludes by showing that $V$ is an $L$-viscosity solution of \eqref{EHJB} for the constant $\lambda$ using the techniques we used in the proof of Proposition \ref{ULip}. In this case there is no time variable but $\lambda^N\to \lambda$. We omit the details.
\vspace{0.5mm}

\nd
\textit{Step 2.} (Uniqueness)\\
    We start by showing the second part. We observe that if $V\colon \mathcal{P}_2\to \R$ is an $L$-viscosity solution of \eqref{EHJB}, then $\Psi(t,m):=V(m)+\lambda t,\; (t,m)\in [0,+\infty)\times\mathcal{P}_2$ satisfies
    \be\label{HJequation}
    \partial_t\Psi(t,m) +\int H(x, D_m\Psi(t,m,x),m)m(dx)=0,\;\;(t,m)\in (0,+\infty)\times\mathcal{P}_2,
    \ee
    in the $L$-viscosity sense. Indeed, suppose that $\Phi\colon (0,+\infty)\times \mathcal{H}\to \R$ be a $C^{1,1}$ test function such that $\hat{\Psi}(t,X)-\Phi(t,X)$ has a maximum at $(t_0,X_0)\in (0,+\infty)\times \mathcal{H}$. Then,
    \begin{align*}
        &\hat{V}(X_0)+\lambda t-\Phi(t,X_0)\text{ has a max at }t_0>0\implies \lambda=\partial_t\Phi(t_0,X_0);\\
        & \hat{V}(X)+\lambda t_0-\Phi(t_0,X)\text{ has a max at }X_0\in\mathcal{H}\implies \lambda+ \hat{H}(X_0,D\Phi(t_0,X_0))\le 0,
    \end{align*}
    by Definition \ref{evisc.}. Combining both relations yields $\partial_t\Phi(t_0,X_0)+\hat{H}(X_0,D\Phi(t_0,X_0))\le 0$, hence $\Psi$ is viscosity subsolution in the sense of Definition \ref{visc.}. The supersolution property can be proved similarly.
    \vspace{1mm}

    \nd
    It follows that $\Psi(t,m)\pm \sup_{m\in \mathcal{P}_2} \left\{|G(m)-V(m)|\right\}$ also satisfy \eqref{HJequation} in the $L$-viscosity sense. Moreover, setting $t=0$, we observe that
    $$\Psi(0,m)-\sup_{m\in \mathcal{P}_2} \left\{|G(m)-V(m)|\right\}\le U(0,m)=G(m)\le \Psi(0,m)+\sup_{m\in \mathcal{P}_2} \left\{|G(m)-V(m)|\right\}.$$
    Thus, by the comparison principle (Theorem \ref{comp})
    $$|U(t,m)-\Psi(t,m)|\le \sup_{m\in \mathcal{P}_2} \left\{|G(m)-V(m)|\right\}.$$
    This means that $U(t,m)-V(m)-\lambda t$ is uniformly bounded, which implies that $U(t,m)-\lambda t$ is uniformly bounded. By dividing by $t>0$ and setting $t\to +\infty$ we obtain that $U(t,m)/t\xrightarrow{t\to +\infty}\lambda$, hence $\lambda$ is uniquely determined.
    \vspace{1mm}

    \nd
    For the first part, let $(V_1,\lambda_1)$ and $(V_2,\lambda_2)$ be two solutions of \eqref{EHJB}. We observe that $\Psi_1(t,m):=V_1(m)+\lambda_1t$ and $\Psi_2(t,m):=V_2(m)+\lambda_2t$, $(t,m)\in [0,+\infty)\times\mathcal{P}_2$ satisfy \eqref{HJequation} in the viscosity sense and
    $$\Psi_2(0,m)-\sup_{m\in\mathcal{P}_2}\left\{|V_1(m)-V_2(m)|\right\}\le \Psi_1(0,m)\le \Psi_2(0,m)+\sup_{m\in\mathcal{P}_2}\left\{|V_1(m)-V_2(m)|\right\}.$$
    The comparison principle gives
    $$|\Psi_1(t,m)-\Psi_2(t,m)|\le \sup_{m\in\mathcal{P}_2}\left\{|V_1(m)-V_2(m)|\right\},\;\;\text{for all }(t,m)\in [0,+\infty)\times\mathcal{P}_2.$$
    Dividing by $t>0$ and sending $t\to +\infty$ gives $|\lambda_1-\lambda_2|\le 0$, that is $\lambda_1=\lambda_2$.
\end{proof}

\begin{remark}
    The above Proposition shows that $U(t,m)/t$ converges to a constant as $t\to +\infty$. We refer to Proposition \ref{erglimit} for a different proof of this fact.
\end{remark}

\subsection{Proof of Theorem \ref{main1}}\label{4.2}
The previous subsection established that $U(t,m)-\lambda t$ is bounded, if $\lambda$ is the unique ergodic constant. In order to establish its asymptotic behavior as $t\to +\infty$, we prove the Asymptotic Monotonicity Property (Proposition \ref{AMP})

\begin{proof}[Proof of Proposition \ref{AMP}]
Throughout the proof the ergodic constant of Proposition \ref{ergodicconstant} is $\lambda=0$, and $V$ denotes a ${\bf d}_2$-Lipschitz $L$-viscosity solution of \eqref{EHJB} with $\lambda=0$. We write $C_0$ for the Lipschitz constant of $U$ (Proposition \ref{ULip}) and $L_V$ for that of $V$. Both are Lipschitz constants with respect to ${\bf d}_2$, so the lifts satisfy
\be\label{liplifts}
|\hat U(t,X)-\hat U(s,Y)|\le C_0\big(|t-s|+\|X-Y\|_2\big),\qquad |\hat V(X)-\hat V(Y)|\le L_V\|X-Y\|_2 .
\ee
Since $U$ and $V$ are bounded, after adding a constant to $V$ we may assume that
\be\label{UVbounds}
1\le U(t,m)-V(m)\le C_*,\qquad\text{for all }(t,m)\in[0,+\infty)\times\mathcal P_2,
\ee
for some constant $C_*\ge1$, and we set $C:=C_*+1$. Fix $\eta>0$ and define, for $s\ge0$,
$$\mu(s):=\inf_{m\in\mathcal P_2,\;t\ge s}\left\{\frac{U(t,m)-V(m)+\eta(t-s)}{U(s,m)-V(m)} \right\}.$$
By \eqref{UVbounds} every quotient is at least $C_*^{-1}$, and the choice $t=s$ gives the value $1$; hence $C_*^{-1}\le\mu(s)\le1$.
\vspace{1mm}

\nd
We observe that $\mu(s)$ is Lipschitz. Indeed, if $t-s>C_*/\eta$, the numerator in the definition of $\mu$ is larger than $1+C_*>C_*$, which exceeds the denominator; so the quotient is larger than $1$. Writing $\tau:=t-s$, we therefore have
\be\label{mu-compact}
\mu(s)=\min_{m\in\mathcal P_2,\;\tau\in[0,C_*/\eta]}q(s,m,\tau),\qquad\text{for }\;\; q(s,m,\tau):=\frac{U(s+\tau,m)-V(m)+\eta\tau}{U(s,m)-V(m)},
\ee
where the minimum is attained because $\mathcal P_2$ is compact and $U,V$ are continuous. For fixed $(m,\tau)$ in this set, the numerator of $q$ lies in $[1,2C_*]$, the denominator in $[1,C_*]$, and both are $C_0$-Lipschitz in $s$ by \eqref{liplifts}. Thus $s\mapsto q(s,m,\tau)$ is $C_0(1+2C_*)$-Lipschitz, uniformly in $(m,\tau)$, and so is $\mu$, being a minimum of such functions.
\vspace{1mm}

\nd
In order to show \eqref{asmonot} we will prove 
\be\label{eq:step1}
\lim_{s\to+\infty}\mu(s)=1 .
\ee
We can then conclude the proof as follows. By the definition of $\mu$, for all $m\in\mathcal P_2$ and $t\ge s\ge0$,
$$U(t,m)-V(m)+\eta(t-s)\ge\mu(s)\big(U(s,m)-V(m)\big),$$
which, by \eqref{UVbounds}, can be rewritten as
$$U(s,m)-U(t,m)+\eta(s-t)\le\big(1-\mu(s)\big)\big(U(s,m)-V(m)\big)\le C_*\big(1-\mu(s)\big).$$
Hence $\delta_\eta(s):=\min\{C_*-1,C_*(1-\mu(s))\}$ has the required properties as soon as \eqref{eq:step1} holds. The rest of the proof is devoted to the proof of \eqref{eq:step1}.
\vspace{2mm}

\nd
Let $\psi:=\psi_{\eta/2}>0$ be the constant from \eqref{SCA} with $\eta/2$ in place of $\eta$, and set $k:=\psi/C$, where $C$ was defined after \eqref{UVbounds}. We claim that it suffices to prove that $\mu$ is a viscosity supersolution of 
\be\label{auxPDE}
\begin{cases}
\max\left\{ w(s)-1, w'(s)+k(w(s)-1)\right\}=0,\;\;\text{in }(0,+\infty),\\
w(0)=\mu(0).
\end{cases}
\ee
If this is true, then since the function $w(s):=1-\big(1-\mu(0)\big)e^{-ks}$ satisfies $w(0)=\mu(0)$ and $w'=k(1-w)$, the comparison principle yields
$$1\ge \mu(s)\ge1-\big(1-\mu(0)\big)e^{-ks}\xrightarrow{s\to+\infty}1 ,$$
which implies \eqref{eq:step1}. We now show the supersolution property.
\vspace{1mm}

\nd
\textit{Step 1.} (Set-up for the viscosity argument)\\
Let $\phi$ and $\bar s$ be such that $\mu-\phi$ has a minimum at $\bar{s}$. If $\mu(\bar s)=1$, then we have nothing to prove, hence we assume $\mu(\bar s)<1$. Replacing $\phi(s)$ by $\phi(s)-|s-\bar s|^2$ plus a constant, which does not change $\phi'(\bar s)$, we may also assume that $\phi(\bar s)=\mu(\bar s)$ and that, for some $r\in(0,\bar s)$,
\be\label{strictmin}
\mu(s)-\phi(s)\ge|s-\bar s|^2,\qquad\text{for }s\in I_r:=[\bar s-r,\bar s+r].
\ee
By \eqref{mu-compact} there are $\bar m\in\mathcal P_2$ and $\bar t\in[\bar s,\bar s+C_*/\eta]$ realizing the minimum defining $\mu(\bar s)$. We let $\bar X\in\mathcal H$ with $\mathcal L(\bar X)=\bar m$.
\vspace{1mm}

\nd
Set $\rho:=\min\{1,(2(C_0+L_V))^{-1}\}$, $T_1:=3C_*/\eta$, and consider the closed set
$$\mathcal D:=\big\{\xi=(X,Y,Z,t,s)\in\mathcal H^3\times\R^2:\ \|X-Y\|_2\le\rho,\ \|X-Z\|_2\le\rho,\ s\in I_r,\ s\le t\le s+T_1\big\}.$$
To simplify the expressions, for $\xi\in\mathcal D$ we put
$$A(\xi):=\hat U(t,X)-\hat V(Z)+\eta(t-s)\quad\text{and}\quad B(\xi):=\hat U(s,Y)-\hat V(Z).$$
Comparing with $\hat U(t,X)-\hat V(X)$ and $\hat U(s,X)-\hat V(X)$ and using \eqref{liplifts}, \eqref{UVbounds} and the choice of $\rho$, we get
\be\label{ABbounds}
\tfrac12\le A(\xi)\le 4C_*+\tfrac12\quad\text{and}\quad \tfrac12\le B(\xi)\le C_*+\tfrac12\le C,\quad\text{for all }\xi\in\mathcal D.
\ee
For $\e\in(0,1]$ and $\kappa:=\e^2$, we define $\Phi\colon\mathcal H^3\times\R^2\to(-\infty,+\infty]$ with
$$\Phi(\xi):=\begin{cases}\frac{A(\xi)}{B(\xi)}-\phi(s)+\frac1{\e^2}\big(\|X-Y\|_2^2+\|X-Z\|_2^2\big),& \text{ if }\xi\in\mathcal D,\\
+\infty,&\text{ otherwise. }\end{cases}$$
$\Phi$ is lower semicontinuous and bounded from below, and $\Phi(\bar X,\bar X,\bar X,\bar t,\bar s)=\mu(\bar s)-\phi(\bar s)=0$, so $\inf\Phi\le0$. By Lemma \ref{BP}, there are $\hat\xi=(\hat X,\hat Y,\hat Z,\hat t,\hat s)\in\mathcal D$ and $\xi'=(X',Y',Z',t',s')$ with $|\hat\xi-\xi'|\le1$ such that $\Phi(\hat\xi)\le\inf\Phi+\kappa\le\kappa$ and
\be\label{minpoint}
\xi\mapsto\Phi(\xi)+\kappa\big(\|X-X'\|_2^2+\|Y-Y'\|_2^2+\|Z-Z'\|_2^2+|t-t'|^2+|s-s'|^2\big)
\ee
attains its minimum at $\hat\xi$. We write
\be\label{not1}
\begin{split}
e_X:=2\kappa(\hat X-&X'),\quad e_Y:=2\kappa(\hat Y-Y'),\quad e_Z:=2\kappa(\hat Z-Z'),\\
& e_t:=2\kappa(\hat t-t'),\quad e_s:=2\kappa(\hat s-s'),
\end{split}
\ee
all of norm at most $2\kappa$, and
\be\label{not2}
\hat a:=A(\hat\xi),\qquad \hat b:=B(\hat\xi)\in[\tfrac12,C],\qquad \hat\mu:=\hat a/\hat b>0,\qquad \hat v:=\hat V(\hat Z).
\ee
In what follows, $K$ denotes a constant depending only on $C_0,L_V,C_*,C_H,\eta$, which may change from line to line; in particular $K$ is independent of $\e$.
\vspace{2mm}

\nd
\textit{Step 2.} ($\hat\xi$ is in the interior of $\mathcal{D}$)\\[0.3ex]
Let $\xi\in\mathcal D$ and $\xi_0:=(X,X,X,t,s)\in\mathcal D$. By \eqref{liplifts}, $|A(\xi)-A(\xi_0)|\le L_V\|X-Z\|_2$ and $|B(\xi)-B(\xi_0)|\le C_0\|X-Y\|_2+L_V\|X-Z\|_2$, so by \eqref{ABbounds} the quotients $A/B$ at $\xi$ and at $\xi_0$ differ by at most $K(\|X-Y\|_2+\|X-Z\|_2)$. Since $t\ge s$, the quotient at $\xi_0$ is at least $\mu(s)$ by the definition of $\mu$. Hence
\be\label{lowerAB}
\frac{A(\xi)}{B(\xi)}\ge\mu(s)-K\big(\|X-Y\|_2+\|X-Z\|_2\big),\quad\text{for all }\xi\in\mathcal D.
\ee
Applying \eqref{lowerAB} at $\hat\xi$, Young's inequality\footnote{$Ka\le\frac{a^2}{4\e^2}+K^2\e^2$, for any $a\ge0$.}, and \eqref{strictmin}, we obtain
$$\e^2\ge\Phi(\hat\xi)\ge\mu(\hat s)-\phi(\hat s)+\frac{\|\hat X-\hat Y\|_2^2+\|\hat X-\hat Z\|_2^2}{2\e^2}-K^2\e^2\ge|\hat s-\bar s|^2+\frac{\|\hat X-\hat Y\|_2^2+\|\hat X-\hat Z\|_2^2}{2\e^2}-K^2\e^2,$$
therefore
\be\label{closeness}
|\hat s-\bar s|\le K\e\quad\text{and}\quad \|\hat X-\hat Y\|_2+\|\hat X-\hat Z\|_2\le K\e^2 .
\ee
Moreover, \eqref{lowerAB} at $\hat\xi$ and $\Phi(\hat\xi)\le\e^2$ give $\mu(\hat s)-K\e^2\le\hat\mu\le\phi(\hat s)+\e^2$. Thus, by \eqref{closeness} and the continuity of $\mu$ and $\phi$, $\hat \mu\to\mu(\bar s)<1$, as $\e\to0$. Next, we observe that
\begin{itemize}
\item If $\hat t=\hat s$, then $\hat a-\hat b=\hat U(\hat s,\hat X)-\hat U(\hat s,\hat Y)\ge-C_0\|\hat X-\hat Y\|_2$, so $\hat\mu=1+(\hat a-\hat b)/\hat b\ge1-K\e^2$, which contradicts $\hat\mu\to\mu(\bar s)<1$ for $\e$ small.
\item If $\hat t=\hat s+T_1$, then $\hat a\ge\eta T_1=3C_*$ and $\hat\mu\ge3C_*/C\ge\frac32$, which again contradicts $\hat{\mu}<1$ for $\e$ small.
\end{itemize}
Consequently, for $\e$ small enough, $\hat\xi$ lies in the interior of $\mathcal D$, $\hat\mu\le\frac{1+\mu(\bar s)}2<1$, and $\hat t>\hat s\ge\bar s-r>0$. In particular, \eqref{minpoint} has a local minimum at $\hat{\xi}$ in the whole space $\mathcal H^3\times\R^2$.
\vspace{2mm}

\nd
\textit{Step 3.} (Three viscosity inequalities)\\[0.3ex]
We freeze all but one group of the 3 groups of variables in \eqref{minpoint}: $(t,X)$, $Z$, and $(s,X)$. This will allow us to write one viscosity inequality for each group.
\vspace{1mm}

\nd
\emph{(i) $U$ at $(\hat t,\hat X)$.} Since $B$ does not depend on $(t,X)$, the map
$$(t,X)\mapsto\frac{\hat U(t,X)-\hat v+\eta(t-\hat s)}{\hat b}+\frac1{\e^2}\big(\|X-\hat Y\|_2^2+\|X-\hat Z\|_2^2\big)+\kappa\big(\|X-X'\|_2^2+|t-t'|^2\big)$$
has a local minimum at $(\hat t,\hat X)$. Multiplying by $\hat b>0$, we see that $\hat U-\psi_1$ has a local minimum at $(\hat t,\hat X)$, where
$$\psi_1(t,X):=-\eta t-\hat b\Big(\frac1{\e^2}\big(\|X-\hat Y\|_2^2+\|X-\hat Z\|_2^2\big)+\kappa\|X-X'\|_2^2+\kappa|t-t'|^2\Big).$$
Since $\hat t>0$, the supersolution property of $U$ gives $\partial_t\psi_1(\hat t,\hat X)+\hat H(\hat X,D\psi_1(\hat t,\hat X))\ge0$. Recalling the notation from \eqref{not1}, \eqref{not2}, this gives
\be\label{ineq1}
\hat H(\hat X,P_t)\ge\eta+\hat b\,e_t\ge\eta-2C\kappa,\;\;\text{where}\quad P_t:=\frac{2\hat b}{\e^2}(\hat Y+\hat Z-2\hat X)-\hat b\,e_X .
\ee

\nd
\emph{(ii) $V$ at $\hat Z$.} Set $a_0:=\hat U(\hat t,\hat X)+\eta(\hat t-\hat s)$ and $b_0:=\hat U(\hat s,\hat Y)$, so that $\hat a=a_0-\hat v$ and $\hat b=b_0-\hat v$, and let $g(Z):=\frac1{\e^2}\|\hat X-Z\|_2^2+\kappa\|Z-Z'\|_2^2$. By the local minimum in \eqref{minpoint}, for $Z$ near $\hat Z$,
$$f\big(\hat V(Z)\big)\ge \hat\mu+g(\hat Z)-g(Z),\qquad\text{where}\quad f(v):=\frac{a_0-v}{b_0-v}=1-\frac{b_0-a_0}{b_0-v}.$$
Here $b_0-a_0=\hat b-\hat a=(1-\hat\mu)\hat b>0$, while $\hat V(Z)<b_0$ and $\hat\mu+g(\hat Z)-g(Z)<1$, for $Z$ near $\hat Z$, due to continuity. For such $Z$, the inequality $f(\hat V(Z))\ge \hat\mu+g(\hat Z)-g(Z)$ is equivalent to
$$\hat V(Z)\le\psi_2(Z):=b_0-\frac{b_0-a_0}{1-(\hat\mu+g(\hat Z)-g(Z))} .$$
Since $\psi_2(\hat Z)=b_0-\hat b=\hat v$, the function $\hat V-\psi_2$ has a local maximum at $\hat Z$, and
\be\label{Q}
D\psi_2(\hat Z)=\frac{b_0-a_0}{(1-\hat\mu)^2}\,Dg(\hat Z)=\frac{\hat b}{1-\hat\mu}\Big(\frac2{\e^2}(\hat Z-\hat X)+e_Z\Big)=:Q .
\ee
Since $\lambda=0$, the subsolution property of $V$ gives
\be\label{ineq2}
\hat H(\hat Z,Q)\le0 .
\ee

\nd
\emph{(iii) $U$ at $(\hat s,\hat Y)$.} Set $a_1:=\hat U(\hat t,\hat X)-\hat v+\eta\hat t$, so that $a_1-\eta\hat s=\hat a$, and let $g_3(s,Y):=\frac1{\e^2}\|\hat X-Y\|_2^2+\kappa\|Y-Y'\|_2^2+\kappa|s-s'|^2$. By the local minimum in \eqref{minpoint}, for $(s,Y)$ near $(\hat s,\hat Y)$,
$$\frac{a_1-\eta s}{\hat U(s,Y)-\hat v}\ge \hat\mu+\phi(s)-\phi(\hat s)+g_3(\hat s,\hat Y)-g_3(s,Y).$$
At $(\hat s,\hat Y)$ the quantities $a_1-\eta s$, $\hat U(s,Y)-\hat v$ and the right hand side of the above inequality equal $\hat a,\hat b,\hat\mu$ respectively, hence they are positive near $(\hat s,\hat Y)$ by continuity. This implies that the inequality is equivalent to 
$$\hat{U}(s,Y)\le\psi_3(s,Y):=\hat v+\frac{a_1-\eta s}{\hat\mu+\phi(s)-\phi(\hat s)+g_3(\hat s,\hat Y)-g_3(s,Y)} .$$
Since $\psi_3(\hat s,\hat Y)=\hat v+\hat a/\hat\mu=\hat U(\hat s,\hat Y)$, the function $\hat U-\psi_3$ has a local maximum at $(\hat s,\hat Y)$. A direct computation gives
\be\label{Ps}
\partial_s\psi_3(\hat s,\hat Y)=-\frac{\eta+\hat b\,\big(\phi'(\hat s)-e_s\big)}{\hat\mu},\qquad D\psi_3(\hat s,\hat Y)=\frac{\hat b}{\hat\mu}\Big(\frac2{\e^2}(\hat Y-\hat X)+e_Y\Big)=:P_s .
\ee
Since $\hat s>0$, the subsolution property of $U$ gives $\partial_s\psi_3(\hat s,\hat Y)+\hat H(\hat Y,P_s)\le0$, that is,
\be\label{ineq3}
\hat\mu\,\hat H(\hat Y,P_s)\le\eta+\hat b\,\phi'(\hat s)-\hat b\,e_s\le\eta+\hat b\,\phi'(\hat s)+2C\kappa .
\ee

\nd
By the Lipschitz continuity of $U,V$ (see equation \eqref{liplifts}) we also have $\|P_t\|_2,\|P_s\|_2\le C_0$ and $\|Q\|_2\le L_V$, since $P_t,P_s,Q$ (defined in \eqref{ineq1}, \eqref{Q}, \eqref{Ps}, respectively) are gradients of $C^1$ test functions evaluated at touching points with $U,V$.
\vspace{2mm}

\nd
\textit{Step 4.} (Conclusion via \eqref{SCA})\\[0.3ex]
Set $z=\hat b\,(e_X+e_Y+e_Z)$. From the definitions of $P_t,P_s,Q$ we have the identity
\be\label{convexid}
\hat\mu P_s+(1-\hat\mu)Q=P_t+z,\;\;\text{ where }  \bigl\|z\bigr\|_2\le6C\kappa .
\ee
Using Assumption (A1)(1), \eqref{closeness}, the gradient bounds, $\hat\mu\le1$ and $\kappa=\e^2$, the inequalities \eqref{ineq1}, \eqref{ineq2}, and \eqref{ineq3} give
\begin{align}
\hat H(\hat X,Q)&\le\hat H(\hat Z,Q)+C_H(1+2L_V)\|\hat X-\hat Z\|_2\le K\e^2,\label{t2}\\
\hat\mu\,\hat H(\hat X,P_s)&\le\hat\mu\,\hat H(\hat Y,P_s)+C_H(1+2C_0)\|\hat X-\hat Y\|_2\le\eta+\hat b\,\phi'(\hat s)+K\e^2,\label{t3}\\
\hat H(\hat X,P_t+z)&\ge\hat H(\hat X,P_t)-C_H(1+2C_0+6C\kappa)\|z\|_2\ge\eta-K\e^2 .\label{t1}
\end{align}
Choosing $\e$ small, the right hand side of \eqref{t1} is at least $\eta/2$, while the right hand side of \eqref{t2} is at most $\eta/4$, hence
% Set $q:=Q$ and $p:=\hat\mu(P_s-Q)$, so that $p+q=P_t+\mathfrak r$ and $p/\hat\mu+q=P_s$ by \eqref{convexid}. For $\e$ small, \eqref{t1} and \eqref{t2} give $\hat H(\hat X,p+q)\ge\eta/2$ and $0\le\theta\le\eta/4$; since also $\hat\mu\in(0,1)$, 
Assumption (i-SCA) with $\eta/2$ in place of $\eta$ and $\theta= K\e^2$ yields
$$\hat\mu\,\hat H(\hat X,P_s)\ge\hat H\bigl(\hat X,P_t+z\bigr)+(\psi-c_*K\e^2)(1-\hat\mu).$$
Combining this with \eqref{t3} and \eqref{t1} yields $\hat b\,\phi'(\hat s)\ge(\psi-c_*K\e^2)(1-\hat\mu)-K\e^2 .$
For $\e$ small, the first term on the right is nonnegative, so dividing by $\hat b\in[\frac12,C]$ we get
$$\phi'(\hat s)\ge\frac{(\psi-c_*K\e^2)(1-\hat\mu)}{C}-K\e^2 .$$
Letting $\e\to0$ and using $\hat s\to\bar s$, $\hat\mu\to\mu(\bar s)$ and $\phi\in C^1$, we obtain $\phi'(\bar s)\ge\frac{\psi}{C}\big(1-\mu(\bar s)\big)=k\big(1-\mu(\bar s)\big)$, which is the desired supersolution inequality. The proof is complete.
\end{proof}

\nd
We can now show Theorem \ref{main1}.

\begin{proof}[Proof of Theorem \ref{main1}]
We once again assume, without loss of generality, that $\lambda=0$.
\vspace{1mm}

\nd
\textit{Step 1.} (Limit along a subsequence of times)\\
By the Arzelà-Ascoli theorem on the compact space $\mathcal P_2$, there exists $T_n\to+\infty$ such that $(U(T_n,\cdot))_{n\in\mathbb N}$ converges uniformly. Let $c_{n,k}:=\sup_{m}|U(T_n,m)-U(T_k,m)|$. The functions $U(T_n+\cdot,\cdot)$ and $U(T_k+\cdot,\cdot)\pm c_{n,k}$ are $L$-viscosity solutions of \eqref{HJB0} in $(0,+\infty)\times\mathcal P_2$, and they are ordered at $t=0$. Hence, by the comparison principle,
$$\sup_{t\ge0,\,m\in\mathcal P_2}|U(T_n+t,m)-U(T_k+t,m)|\le c_{n,k}\xrightarrow{n,k\to+\infty}0 .$$
Thus $U(T_n+\cdot,\cdot)$ converges uniformly on $[0,+\infty)\times\mathcal P_2$ to some $U^\infty$, which satisfies $|U^\infty|\le C$, for some constant $C>0$, and is $C_0$-Lipschitz in $m$.
\vspace{2mm}

\nd
\textit{Step 2.} ($U^\infty$ is nondecreasing in time)\\
Let $m\in\mathcal P_2$, $t\ge s\ge0$ and $\eta>0$. By \eqref{asmonot} at the times $s+T_n\le t+T_n$,
$$U(s+T_n,m)-U(t+T_n,m)+\eta(s-t)\le\delta_\eta(s+T_n).$$
Letting $n\to+\infty$ and then $\eta\to0$, we obtain $U^\infty(s,m)\le U^\infty(t,m)$. Hence $U^\infty(t,m)$ increases to a limit $V(m)$ as $t\to+\infty$, and $V$ is $C_0$-Lipschitz. Since $\mathcal P_2$ is compact and $V$ is continuous, Dini's theorem shows that the convergence is uniform in $m$.
\vspace{2mm}

\nd
\textit{Step 3.} (Convergence of $U$ and viscosity properties of the limit)\\
Let $\e>0$. Choose $n$ such that $\sup_{t,m}|U(T_n+t,m)-U^\infty(t,m)|\le\e$, and then $\tau$ such that $\sup_m|U^\infty(t,m)-V(m)|\le\e$ for $t\ge\tau$. Then $|U(t,m)-V(m)|\le2\e$ for all $t\ge T_n+\tau$ and $m\in\mathcal P_2$. Hence \eqref{LTB} holds with $\lambda=0$. It remains to show that $V$ solves \eqref{EHJB} with $\lambda=0$.
\vspace{2mm}

\nd
For any $k\in \mathbb{N}$, it follows from the previous paragraph that $e_k:=\sup_{t\ge k,\,m\in\mathcal P_2}|U(t,m)-V(m)|\xrightarrow{k\to \infty}0$. Let $\varphi$ be a test function such that $\hat V-\varphi$ has a maximum at $X_0$ on $\bar B(X_0,r)$, for some $r\in(0,1)$. Set $z_k:=e_k+1/k$, and define,
$$\Phi_k(t,X):=\begin{cases}-\hat U(t+k,X)+\varphi(X)+|t-1|^2+\|X-X_0\|_2^2,&\text{on }[1-r,1+r]\times\bar B(X_0,r),\\
+\infty,&\text{otherwise.}
\end{cases}$$
By Lemma \ref{BP}, there exist $(t_k,X_k)$ and $(t_k',X_k')$ so that $\Phi_k+z_k\big(|t-t_k'|^2+\|X-X_k'\|_2^2\big)$ attains a minimum at $(t_k,X_k)$, $\Phi_k(t_k,X_k)\le\Phi_k(1,X_0)+z_k$ and
$|2z_k(t_k-t_k')|,\|2z_k(X_k-X_k')\|_2\le2z_k$.
\vspace{1mm}

\nd
We have $|\hat U(t+k,\cdot)-\hat V|\le e_k$ and $(\hat V-\varphi)(X_k)\le(\hat V-\varphi)(X_0)$. Hence the first property gives
$$|t_k-1|^2+\|X_k-X_0\|_2^2\le2e_k+z_k\xrightarrow{k\to+\infty}0 .$$
So, for $k$ large, $(t_k,X_k)$ is an interior point, and 
$$(t,X)\mapsto\hat U(t+k,X)-\left(\varphi(X)+|t-1|^2+\|X-X_0\|_2^2+z_k\big(|t-t_k'|^2+\|X-X_k'\|_2^2\big)\right)$$
has a local maximum at $(t_k,X_k)$. The subsolution property of $\hat U(\cdot+k,\cdot)$ gives
$$2(t_k-1)+2z_k(t_k-t_k')+\hat H\big(X_k,\,D\varphi(X_k)+2(X_k-X_0)+2z_k(X_k-X_k')\big)\le0 .$$
Letting $k\to+\infty$ and using Assumption (A1)(1), we obtain $\hat H(X_0,D\varphi(X_0))\le0$. The supersolution property is proved in the same way, applying Lemma \ref{BP} to $\hat U(t+k,X)-\varphi(X)+|t-1|^2+\|X-X_0\|_2^2$.
\end{proof}

\section{Counterexample in $\R^d$}\label{counter}

\nd
If \eqref{HJB0} is posed in $\mathcal{P}_2(\R^d)$ (instead of $\mathcal{P}_2(\T^d)$), so $H\colon \R^d\times \R^d\times\mathcal{P}_2(\R^d)\to \R$ and $G\colon \mathcal{P}_2(\R^d)\to \R$ satisfying Assumptions (A1)-(A2) with $\T^d$ replaced by $\R^d$, then our main result (Theorem \ref{main1}) does not hold even if $p\mapsto H(x,p,m)$ is strictly convex. We present the following example.

\begin{example}
Fix $b\in\R^d$ with $|b|=1$, and set
\be\label{data}\nonumber
H(x,p,m):=\frac12|p|^2+b\cdot p,\qquad G(m):=\int_{\R^d}g(y)\,m(dy)\;\;\text{with }\;\; g(y):=\sin\big(\log(1+|y|)\big).
\ee
Assumptions (A1)-(A2) are clearly satisfied. By Lemma \ref{lem:linearG} below, the unique $L$-viscosity solution of \eqref{HJB0} is given by
$$U(t,m)=\int_{\R^d} u(t,x)\;m(dx),\;\;\text{where}\;\;u(t,x)=\min_{y\in\R^d}\left\{g(y)+\frac{|y-(x-bt)|^2}{2t}\right\}$$
is the solution of the finite dimensional HJ equation $\partial_tu+H(x,Du)=0,$ for $(t,x)\in (0,+\infty)\times \R^d$, and $u(0,x)=g(x)$, for $x\in \R^d$ and is obtained by the Hopf-Lax formula \cite[Sections 3.3 and 10.3]{evans2022partial}. It is straightforward to check that for fixed $x\in \R^d$ and $t$ large enough
$$g(x-bt)-\frac{4}{\sqrt t}\le u(t,x)\le g(x-bt).$$
Since $|g(x-bt)-g(-bt)|\le\frac{2|x|}{t}$, for $t\ge 2|x|$, we may remove the $x$-dependence of $g$ in the above inequality and obtain $|u(t,x)-g(-bt)|\xrightarrow{t\to+\infty}0$, for any $x\in \R^d$. This directly implies
$$\lim_{t\to+\infty}\Big(U(t,m)-\sin\big(\log(1+t)\big)\Big)=0,$$
hence $\liminf_{t\to+\infty}U(t,m)=-1$ and $\limsup_{t\to+\infty}U(t,m)=1$, so $U(t,m)-\lambda t$ does not converge to a function $V\colon\mathcal{P}_2(\R^d)\to \R$ as $t\to+\infty$, for any $m\in\mathcal{P}_2(\R^d)$ and $\lambda\in \R$.
\end{example}

%\textcolor{cyan}{\textbf{Discussion:} Do we want to include more details? Reasons that this happens? This will require us to say more about the connection with control problems (in the strongly convex case).}
\begin{remark}[Outlook]
In the above example, for $T>0$, the function $W^T(t,m):=U(T-t,m)$ is the value function of a deterministic optimal control problem of the form \eqref{MFCvalue} with reaction cost $L(x,a,m)=\frac{1}{2}|a+b|^2$ and terminal cost $G(m)=\int g(x) m(dx)$. This shows that additional assumptions (confinement mechanisms) are therefore needed in order for long time behavior to appear in problems posed in $\R^d$. We expect that such examples include strictly convex running $L$ and terminal $G$ costs or a confining drift (the latter being the case in \cite{fuhrman2025ergodic}). We intend to study such extensions, including mean field control problems (recall \eqref{MFCvalue}) with stochastic dynamics, in future works. 
\end{remark}
%\vspace{1mm}

\begin{lemma}\label{lem:linearG}
Let $b\in\R^d$, let $g\colon\R^d\to\R$ be bounded and Lipschitz continuous, and let $H(x,p,m):=\frac12|p|^2+b\cdot p$ and $G(m):=\int_{\R^d}g\,dm$. Let $u\colon [0,+\infty)\times \R^d\to \R$ be the unique bounded viscosity solution of
$$\partial_t u+H(x,Du)=0,\;\;(t,x)\in (0,+\infty)\times \R^d,\;\;\;\; u(0,x)=g(x),\;\; x\in \R^d.$$
Set $U(t,m):=\int_{\R^d}u(t,x)\,m(dx)$, $(t,m)\in [0,+\infty)\times \mathcal{P}_2(\R^d)$. Then, $U$ is the unique bounded $L$-viscosity solution of \eqref{HJB0} on $[0,+\infty)\times\mathcal{P}_2(\R^d)$, where $L$-viscosity solutions for this problem can be understood as in Definition \ref{visc.}.
\end{lemma}
 
\begin{proof}
Fix $N\in\mathbb{N}$ and set $V^N(t,\bm x):=\frac1N\sum_{i=1}^Nu(t,x^i)$ for $\bm x\in(\R^d)^N$. We claim that $V^N$ is the unique viscosity solution of the $\R^d$-analogue of \eqref{NHJB}, namely
$$\partial_tV^N+\frac1N\sum_{i=1}^NH\big(x^i,ND_{x^i}V^N,m^N_{\bm x}\big)=0\;\;\text{in }(0,+\infty)\times(\R^d)^N,\quad V^N(0,\bm x)=G(m^N_{\bm x})=\frac1N\sum_{i=1}^Ng(x^i).$$
The Hamiltonian $\mathcal{H}_N(\bm p):=\frac1N\sum_{i=1}^N\big(\frac12|Np^i|^2+b\cdot Np^i\big)$, where $\bm p=(p^1,\ldots,p^N),$ of this equation is convex and superlinear. Its Legendre transform is
$$\mathcal{L}_N(\bm v)=\sup_{\bm p\in(\R^d)^N}\left\{\sum_{i=1}^Np^i\cdot v^i-\mathcal{H}_N(\bm p)\right\}=\frac1N\sum_{i=1}^N\sup_{q\in\R^d}\Big\{q\cdot v^i-\frac12|q|^2-b\cdot q\Big\}=\frac1N\sum_{i=1}^N\frac12|v^i-b|^2.$$
Hence, by the Hopf-Lax formula \cite[Sections 3.3 and 10.3]{evans2022partial}, the unique viscosity solution is
$$\min_{\bm y\in(\R^d)^N}\Big\{\frac1N\sum_{i=1}^Ng(y^i)+t\,\mathcal{L}_N\Big(\frac{\bm x-\bm y}{t}\Big)\Big\}=\frac1N\sum_{i=1}^N\min_{y^i\in\R^d}\Big\{g(y^i)+\frac{|y^i-(x^i-bt)|^2}{2t}\Big\}=V^N(t,\bm x),$$
where in the last step we used the Hopf-Lax formula for $u$.
\vspace{1mm}

\nd
By the previous observation, $U(t,m^N_{\bm x})=V^N(t,\bm x)$, for every $N\in\mathbb{N}$ and $(t,\bm x)\in[0,+\infty)\times(\R^d)^N$. Moreover, $U(0,\cdot)=G(\cdot)$, and $U$ is ${\bf d}_2$-Lipschitz continuous, because $u$ is Lipschitz continuous on $[0,+\infty)\times\R^d$ (see \cite[Section 3.3]{evans2022partial}), and integrating against an optimal coupling of $m,n$ gives
$$|U(t,m)-U(s,n)|\le \mathrm{Lip}(u)\big(|t-s|+{\bf d}_2(m,n)\big).$$
We now repeat the proof of the viscosity property in Proposition~\ref{ULip}, with $\bar U^N$ replaced by $U$. Since $U$ coincides with $V^N$ on empirical measures, the error term there vanishes, i.e., $e_N=0$. The only compactness used in that argument is that of the sets $K_N=\{(t,\bm x):|t-t_0|\le r,\ \|X_{\bm x}-X_0\|_2\le r\}$, which remain compact in $[0,+\infty)\times(\R^d)^N$, since $\|X_{\bm x}\|_2^2=\frac1N\sum_{i=1}^N|x^i|^2$. Hence $U$ is an $L$-viscosity solution of \eqref{HJB0}. Uniqueness follows from the comparison principle \cite[Theorem 7.4]{gangbo2021finite}.
\end{proof}

\appendix

\section{Technical results}\label{appendix}

\subsection{Variational principle}
\nd
Throughout the paper we use the following variational principle that allows us to obtain minimizers of perturbations of lower semi-continuous and bounded from below functions defined on Hilbert spaces. The Lemma below follows from \cite[Theorem 2.5.3]{borwein2005techniques}.

\begin{lemma}[Borwein--Preiss variational principle]\label{BP}
Let $\mathcal{H}$ be a Hilbert space with norm $\|\cdot\|$ and suppose that $\Phi\colon\mathcal{H}\to(-\infty,+\infty]$ is proper
\footnote{Proper in this context means that $\Phi$ is not identically $+\infty$.},
lower semicontinuous and bounded from below. Then for every $\kappa>0$, there exist
$\hat\xi,\xi'\in\mathcal{H}$ with $\|\hat\xi-\xi'\|\le1$ such that $\Phi(\hat\xi)\le\inf\Phi+\kappa$ and
\[
\Phi(\xi)+\kappa\|\xi-\xi'\|^2>\Phi(\hat\xi)+\kappa\|\hat\xi-\xi'\|^2,
\qquad\text{for all }\xi\neq\hat\xi .
\]
\end{lemma}

\subsection{Technical proofs}

\begin{proposition}\label{erglimit}
    Suppose that Assumptions (A1)-(A2) are in place and let $U$ be the $L$-viscosity solution of \eqref{HJB0}. There exists a constant $c\in \R$ such that 
    \be\label{ergconstant}
    \frac{U(t,m)}{t}\xrightarrow{t\to +\infty}c,\;\;\text{ uniformly with respect to }m\in \mathcal{P}_2.
    \ee
\end{proposition}

\begin{proof}
    We set $M(t):=\sup_{m\in \mathcal{P}_2}\{ U(t,m)-G(m)\}$. We note that $M(t)<\infty$ for each $t>0$, due to Proposition \ref{ULip} and the boundedness of $G$. In addition, the Lipschitz continuity of $U$ gives $M(t)\ge -C_0t$.
    \vspace{1mm}

    \nd
    For $t,s>0$, we observe that the functions $(t,m)\mapsto U(t+s,m)$ and $(t,m)\mapsto U(t,m)+M(s)$ are $L$-viscosity solutions of \eqref{HJB0} with initial conditions $U(s,m)\le U(0,m)+M(s)$. Thus, by the comparison principle $U(t+s,m)\le U(t,m)+M(s)$, for all $(t,m)\in [0,+\infty)\times \mathcal{P}_2$ and hence $$\sup_{m\in \mathcal{P}_2}\{ U(t+s,m)-U(t,m)\}\le M(s).$$
    Using this inequality, we get
    \begin{align*}
        M(t+s)&\le \sup_{m\in \mathcal{P}_2}\{ U(t+s,m)-U(t,m)\}+\sup_{m\in\mathcal{P}_2}\{U(t,m)-G(m)\}\le M(s)+M(t).
    \end{align*}
    Since $M(t)$ is subadditive and $M(t)/t$ is lower bounded, Fekete's lemma yields
    $$\frac{M(t)}{t}\xrightarrow{t\to+\infty}c:=\inf_{s>0}\left\{\frac{M(s)}{s}\right\}\in \R.$$
    We finish the proof of \eqref{ergconstant} by noting that for any $(t,m)\in [0,+\infty)\times\mathcal{P}_2$ we have $U(t,m)-M(t)\le G(m)\le \|G\|_\infty$ and
    $$U(t,m)-M(t)\ge U(t,m)-1-\left(U(t,m_1)-G(m_1)\right)\ge -\|G\|_\infty-1-C_0{\bf d}_2(m,m_1),$$
    for some $m_1\in\mathcal{P}_2$ depending on $t$. Thus, $U(t,m)-M(t)$ is upper and lower bounded uniformly in $m$. This gives
    $$\frac{U(t,m)}{t}=\frac{M(t)}{t}+\frac{U(t,m)-M(t)}{t}\xrightarrow{t\to+\infty}c+0=c.$$
\end{proof}

\begin{proof}[Proof of Proposition \ref{VNLip}]
We argue as follows. Since $G$ is ${\bf d}_2$-Lipschitz, by  \cite[Theorem 4.2]{cosso2023smooth}\footnote{Even though functions are defined in $\mathcal{P}_2(\R^d)$ in this reference, the result holds with functions defined in $\mathcal{P}_2(\T^d)$.}, we may find a sequence of functions $G_{k}\colon\mathcal{P}_2\to \R$, $k\in\mathbb{N}$, which are smooth and ${\bf d}_2$-Lipschitz continuous with the same constant as $G$ such that
$$\lim_{k\to \infty}\sup_{m\in B}|G_k(m)-G(m)|=0,$$
for every compact set $B\subset\mathcal{P}_2$. We will show the desired estimate when the initial condition in \eqref{NHJB} is $V^N_{k}(0,\bm x)=G_k(m^N_{\bm x})$ for any $k\in\mathbb{N}$ and then the result will follow by sending $k\to \infty$ due to the stability of viscosity solutions.
%To this end, we first observe that the Lipschitz continuity of $G_k$ gives
%\begin{align*}
%    \left|\frac{G_k(m^N_{\bm x+he^i})-G_k(m^N_{\bm x}) }{|h|}   \right|\le \frac{C_G}{|h|}{\bf d}_2(m^N_{\bm x+he^i},m^N_{\bm x})\le \frac{C_G}{|h|}\sqrt{\frac{|h|^2}{N}}\le C_G.
%\end{align*}
%By sending $|h|\to 0$, we get
Let $m\in \mathcal{P}_2$, $\phi:\T^d\to \R^d$ measurable so that $\phi\in L^2_m$, and set $m_t:=(\text{Id}+t\phi)_\#m$, for $t>0$. We have
\begin{align*}
    \left|\frac{G_k(m_t)-G_k(m)}{t}\right|\le \frac{C_G}{t}{\bf d}_2(m_t,m)\le \frac{C_G}{t}\|t\phi\|_{L^2_m}=C_G\|\phi\|_{L^2_m}.
\end{align*}
Sending $t\to 0^+$ yields
$$\left|\int_{\T^d} D_mG_k(m,x)\cdot \phi(x) m(dx)\right|\le C_G\|\phi\|_{L^2_m}\stackrel{\phi\text{ arbitrary}}{\implies} \int_{\T^d}|D_mG_k(m,x)|^2m(dx)\le C_G^2.$$
Therefore, if $m=m^N_{\bm x}$ for some $\bm x\in (\T^d)^N$, we derive
$$\frac{1}{N}\sum_{i=1}^N|D_mG_k(m^N_{\bm x},x^i)|^2\le C_G^2.$$

\nd
We show the result for the unique solution $V_{k}^N\colon [0,+\infty)\times (\T^d)^N\to \R$ of the PDE
\be\label{HJBNrevtime}
\begin{cases}
    \partial_tV^N_{k}+\frac{1}{N}\sum_{i=1}^NH(x^i,ND_{x^i}V^N_{k},m_{\bm x}^N)=0, &(t,\bm x)\in (0,+\infty)\times (\T^d)^N,\\
    V^N_{k}(0,\bm x)= G_k(m_{\bm x}^N), & \bm x\in (\T^d)^N.
\end{cases}
\ee
Note that, in view of $H$ having quadratic growth, for all $\bm x\in (\T^d)^N$ we have
\begin{align*}
    %c_H=\frac{Nc_H(1+|0|^2)}{N}\le
    \left|\frac{1}{N}\sum_{i=1}^NH\left(x^i,D_mG_k(m^N_{\bm x},x^i),m_{\bm x}^N\right)\right|&\le \frac{N\max\{ K_0, k_0\}+C_H\sum_{i=1}^N\left|D_mG_k(m^N_{\bm x},x^i)\right|^2}{N}\\
    &\le \max\{ K_0, k_0\}+C_HC_G^2.
\end{align*}
This implies $V^{N,+}_{k}(t,\bm x)=Ct+G_k(m^N_{\bm x})$ and $V^{N,-}_{k}(t,\bm x)=-Ct+G_k(m^N_{\bm x})$ are a classical supersolution and a classical subsolution of \eqref{HJBNrevtime}, for some constant $C>0$ large enough and independent of $N,k$, respectively. By comparison principle, we get
$$-Ct+G_k(m^N_{\bm x})\le V^N_{k}(t,\bm x)\le Ct+G_k(m^N_{\bm x}),\;\text{for all }(t,\bm x)\in (0,+\infty)\times (\T^d)^N,$$
which gives the Lipschitz regularity in time at $t=0$. To get the Lipschitz regularity in time for any $s>0$, we notice that $(t,\bm x)\mapsto V^N_{k}(s+t,\bm x)$ is also a viscosity solution of \eqref{HJBNrevtime} with initial condition $V^{N}_{k}(s,\bm x)$. Then, in light of the inequality
$$V^N_{k}(s,\bm x)-\|G_k(m^N_{\cdot})-V^N_{k}(s,\cdot)\|_{\infty}\le G_k(m^N_{\bm x})\le V^N_{k}(s,\bm x)+\|G_k(m^N_{\cdot})-V^N_{k}(s,\cdot)\|_{\infty}$$
and the comparison principle, we get
$$V^N_{k}(t+s,\bm x)-\|G_k(m^N_{\cdot})-V^N_{k}(s,\cdot)\|_{\infty}\le V^N_{k}(t,\bm x)\le V^N_{k}(t+s,\bm x)+\|G_k(m^N_{\cdot})-V^N_{k}(s,\cdot)\|_{\infty}$$
or
$$\left|V^N_{k}(t+s,\bm x)- V^N_{k}(t,\bm x)   \right|\le \|G_k(m^N_{\cdot})-V^N_{k}(s,\cdot)\|_{\infty}\le Cs,$$
for any $(t,\bm x)\in [0,+\infty)\times (\T^d)^N$.
\vspace{2mm}

\nd
We further know that $V^N_{k}$ is Lipschitz continuous in both variables (see \cite[Theorem 1.34]{tran2021hamilton}), therefore $V^N_{k}$ satisfies \eqref{HJBNrevtime} a.e. Thus, we have from the Lipschitz continuity with respect to $t$
\begin{align*}
    C\ge |\partial_tV^N_{k}(t,\bm x)|\ge \frac{1}{N}\sum_{i=1}^NH\left(x^i,ND_{x^i}V^N_{k},m_{\bm x}^N\right)&\ge \frac{1}{N}\left(-Nk_0+c_H\sum_{i=1}^N|ND_{x^i}V^N_{k}(t,\bm x)|^{2}\right)\\
    &\ge -k_0 +c_HN\sum_{i=1}^N|D_{x^i}V^N_{k}(t,\bm x)|^{2},\;\;\text{ a.e.}
\end{align*}
This gives a gradient bound of the form $C_0\ge N\sum_{i=1}^N|D_{x^i}V^N_{k}(t,\bm x)|^{2}$ for almost every $(t,\bm x)\in [0,+\infty)\times (\T^d)^N$. By Fubini's theorem, for a.e. $t\ge 0$, this estimate holds for a.e. $\bm x\in (\T^d)^N$. Let $\rho\in C_c^\infty((\R^d)^N)$ be a standard mollifier supported in the unit ball, set $\rho_\e(\bm z):=\e^{-dN}\rho(\bm z/\e)$ and $V^N_{k,\e}(t,\cdot):=V^N_k(t,\cdot)*\rho_\e$. Since $V^N_k(t,\cdot)$ is Lipschitz, its a.e.\ gradient coincides with its distributional gradient, so that $D_{x^i}V^N_{k,\e}(t,\cdot)=D_{x^i}V^N_k(t,\cdot)*\rho_\e$, $i\in\{1,\ldots,N\}$. Hence, by Minkowski's inequality and the gradient bound,
\be\label{gradbdeps}
\left(\sum_{i=1}^N|D_{x^i}V^N_{k,\e}(t,\bm z)|^2\right)^{\frac12}\le\int_{(\R^d)^N}\left(\sum_{i=1}^N|D_{x^i}V^N_k(t,\bm z-\bm y)|^2\right)^{\frac12}\rho_\e(\bm y)\,d\bm y\le \frac{C_0^{\frac12}}{\sqrt{N}},\; \text{for all }\bm z\in(\R^d)^N.
\ee
Now let $\bm x,\bm y\in(\T^d)^N$. Since the optimal coupling between $m^N_{\bm x}$ and $m^N_{\bm y}$ is induced by a permutation and $V^N_k$ is symmetric in its variables, we may assume that ${\bf d}_{2}^2(m^N_{\bm x}, m^N_{\bm y})=N^{-1}\sum_{i=1}^N|x^i-y^i|_{\T^d}^2$. We choose representatives $\tilde{\bm x},\tilde{\bm y}\in(\R^d)^N$ of $\bm x,\bm y$ such that $|\tilde x^i-\tilde y^i|=|x^i-y^i|_{\T^d}$, so that $\|\tilde{\bm x}-\tilde{\bm y}\|_2=N^{\frac12}{\bf d}_2(m^N_{\bm x},m^N_{\bm y})$, where $\|\cdot\|_2$ denotes the Euclidean norm on $(\R^d)^N$. Since $V^N_{k,\e}(t,\cdot)$ is smooth, \eqref{gradbdeps} gives
\begin{align*}
    |V^N_{k,\e}(t,\tilde{\bm x})-V^N_{k,\e}(t,\tilde{\bm y})|&= \Big|\int_0^1\frac{\text{d}}{\text{d}s} V^N_{k,\e}\big(t, s\tilde{\bm x}+(1-s)\tilde{\bm y}\big)\,ds \Big|\le \sup_{\bm z\in(\R^d)^N}\;\sum_{i=1}^N|D_{x^i}V^N_{k,\e}(t,\bm z)|\;|\tilde x^i-\tilde y^i|\\
    &\le \sup_{\bm z\in(\R^d)^N}\;\left(\sum_{i=1}^N |D_{x^i}V^N_{k,\e}(t,\bm z)|^{2}  \right)^{\frac{1}{2}}\|\tilde{\bm x}-\tilde{\bm y}\|_2\\
    &\le C_0^{\frac{1}{2}}N^{-\frac{1}{2}}\|\tilde{\bm x}-\tilde{\bm y}\|_2=C_0^{\frac{1}{2}}{\bf d}_2(m^N_{\bm x},m^N_{\bm y}).
\end{align*}
Since $V^N_k(t,\cdot)$ is Lipschitz, $|V^N_{k,\e}(t,\cdot)-V^N_k(t,\cdot)|\le \text{Lip}\big(V^N_k(t,\cdot)\big)\,\e$ on $(\R^d)^N$, hence $V^N_{k,\e}(t,\cdot)\to V^N_k(t,\cdot)$ uniformly as $\e\to0$. Letting $\e\to0$ in the above inequality and using the periodicity of $V^N_k$, we obtain
$$|V^N_{k}(t,\bm x)-V^N_{k}(t,\bm y)|\le C_0^{\frac{1}{2}}{\bf d}_2(m^N_{\bm x},m^N_{\bm y}),\quad \bm x,\bm y\in(\T^d)^N,\;\;\text{for a.e. }t\in [0,T].$$
We conclude by the continuity of $V^N_k$ in $t$.
\end{proof}

% Now let $\bm x,\bm y\in (\T^d)^N$ and $t\ge 0$. Since $V^N_{k}$ is symmetrical in the $\bm x$ variables, we may assume that ${\bf d}_{2}^2(m^N_{\bm x}, m^N_{\bm y})=N^{-1}\sum_{i=1}^N|x^i-y^i|^2=N^{-1}\|\bm x-\bm y\|_2^2$, where the distances $|x^i-y^i|,\; i\in \{1,\ldots,N\}$ are understood as distances in $\T^d$. We have
% \begin{align*}
%     |V^N_{k}(t,\bm x)-V^N_{k}(t,\bm y)|= \Big|\int_0^1\frac{\text{d}}{\text{d}s} V^N_{k}(t, s\bm x+(1-s)&\bm y)ds \Big|\le \esssup_{\bm z\in(\R^d)^N}\;\sum_{i=1}^N|D_{x^i}V^N(t,\bm z)|\;|x^i-y^i|\\
%     &\le \esssup_{\bm z\in(\R^d)^N}\;\left(\sum_{i=1}^N |D_{x^i}V^N(t,\bm z)|^{2}  \right)^{\frac{1}{2}}\|\bm x-\bm y\|_2\\
%     &\le C_0^{\frac{1}{2}}N^{-\frac{1}{2}}\|\bm x-\bm y\|_2=C_0^{\frac{1}{2}}{\bf d}_2(m^N_{\bm x},m^N_{\bm y}).
% \end{align*}
% The proof is complete.

\vspace{1mm}

\nd
We finally prove Proposition \ref{ULip} which establishes the existence of an $L$-viscosity solution of \eqref{HJB0}. The proof uses the ideas developed in \cite[Theorem 2.1]{gangbo2021finite} adapted to our setup.

\begin{proof}[Proof of Proposition \ref{ULip}]
We denote by $\mathcal P^N\subset\mathcal P_2$ the set of empirical measures $m^N_{\bm x}$, $\bm x\in(\T^d)^N$. 
\vspace{1mm}
 
\nd
By Proposition \ref{VNLip}, the function $U^N(t,m^N_{\bm x}):=V^N(t,\bm x)$ is well defined and $C_0$-Lipschitz on $[0,+\infty)\times\mathcal P^N$. Its extension
$$\bar U^N(t,m):=\inf_{n\in\mathcal P^N}\big\{U^N(t,n)+C_0\,{\bf d}_2(m,n)\big\}$$
has the following three properties: it coincides with $U^N$ in $[0,+\infty)\times\mathcal{P}^N$, it satisfies the Lipschitz estimate of the statement, and it is bounded on $[0,T]\times\mathcal P_2$ uniformly in $N$, for every $T>0$. By the Arzelà--Ascoli theorem and a diagonal argument, we get that $\bar U^N$ converges to some $U$ uniformly on $[0,T]\times\mathcal P_2$ (up to a subsequence that we do not relabel), for every $T>0$. The limit $U$ satisfies the Lipschitz estimate of the statement.
\vspace{1mm}
 
\nd
Let $m\in\mathcal P_2$. Since $\bigcup_N\mathcal P^N$ is dense in $(\mathcal P_2,{\bf d}_2)$ (see \cite[Theorem 6.18]{villani2009optimal}), we can pick $n_N\in\mathcal P^N$ with ${\bf d}_2(m,n_N)\to0$. Then, by the Lipschitz continuity,
$$|\bar U^N(0,m)-G(m)|\le(C_0+C_G)\,{\bf d}_2(m,n_N)\xrightarrow{N\to +\infty} 0,$$
hence $U(0,\cdot)=G$.
\vspace{2mm}

\nd
We now show that $U$ is an $L$-viscosity solution of \eqref{HJB0}. We only verify the subsolution property; the supersolution property is analogous. Let $\Phi\in C^{1,1}([0,+\infty)\times\mathcal H)$ be such that $\hat U-\Phi$ has a local maximum at $(t_0,X_0)\in(0,+\infty)\times\mathcal H$. Adding $|t-t_0|^2+\|X-X_0\|_2^2$ to $\Phi$ does not change its derivatives at $(t_0,X_0)$, so we may assume that, for some $r\in(0,t_0)$,
\be\label{strict}
(\hat U-\Phi)(t,X)\le(\hat U-\Phi)(t_0,X_0)-|t-t_0|^2-\|X-X_0\|_2^2\quad\text{for }|t-t_0|\le r,\ \|X-X_0\|_2\le r .
\ee
Set $Y_0:=D\Phi(t_0,X_0)$. Standard measure theory arguments imply that, since $(\Omega,\mathcal{F},\mathbb{P})$ is atomless, for any $N\in \mathbb{N}$, there exist a partition $A_1,\dots,A_N$ of $\Omega$ with $\mathbb P(A_i)=1/N$, $\bm x_N\in(\T^d)^N$ and $\bm y_N\in(\R^d)^N$ such that
\be\label{discr}
X_{\bm x_N}:=\sum_{i=1}^Nx_N^i\mathbf 1_{A_i}\to X_0\quad\text{and}\quad Y_N:=\sum_{i=1}^Ny_N^i\mathbf 1_{A_i}\to Y_0\quad\text{in }L^2 .
\ee
Let $\mathcal F_N:=\sigma(A_1,\dots,A_N)$ be the $\sigma$-algebra generated by the sets $(A_i)_{i=1}^N$ and for $\bm x\in(\T^d)^N$ consider $X_{\bm x}:=\sum_{i=1}^Nx^i\mathbf 1_{A_i}$. Since $\mathcal L(X_{\bm x})=m^N_{\bm x}$, we have
\be\label{VNU}
V^N(t,\bm x)=\bar U^N(t,m^N_{\bm x})\quad\text{and}\quad\hat U(t,X_{\bm x})=U(t,m^N_{\bm x}).
\ee
We consider the test functions $\varphi^N(t,\bm x):=\Phi(t,X_{\bm x})$, $(t,\bm x)\in [0,+\infty)\times (\T^d)^N$. Since $\bm x\mapsto X_{\bm x}$ is linear and $\partial_{x^i}X_{\bm x}=\mathbf 1_{A_i}$, the chain rule gives $\varphi^N\in C^1$, with 
$$\partial_t\varphi^N(t,\bm x)=\partial_t\Phi(t,X_{\bm x})\quad\text{and}\quad ND_{x^i}\varphi^N(t,\bm x)=N\,\mathbb E\big[D\Phi(t,X_{\bm x})\mathbf 1_{A_i}\big].$$
We notice that the last expression is the value of $\mathbb E[D\Phi(t,X_{\bm x})\,|\,\mathcal F_N]$ on $A_i$. Hence, by the definition of $\hat H$,
\be\label{Hdisc}
\frac1N\sum_{i=1}^NH\big(x^i,ND_{x^i}\varphi^N(t,\bm x),m^N_{\bm x}\big)=\hat H\Big(X_{\bm x},\,\mathbb E\big[D\Phi(t,X_{\bm x})\,|\,\mathcal F_N\big]\Big).
\ee
The set $K_N:=\{(t,\bm x):\,|t-t_0|\le r,\ \|X_{\bm x}-X_0\|_2\le r\}$ is compact, and it contains $(t_0,\bm x_N)$ for $N$ large. Let $(s_N,\bm z_N)$ be a maximum point of $V^N-\varphi^N$ on $K_N$, and set $Z_N:=X_{\bm z_N}$. By \eqref{VNU}, $|V^N(t,\bm x)-\hat U(t,X_{\bm x})|\le e_N:=\sup_{[0,t_0+r]\times\mathcal P_2}|\bar U^N-U|$, on $K_N$. Comparing the values of $V^N-\varphi^N$ at $(s_N,\bm z_N)$ and $(t_0,\bm x_N)$, and using \eqref{strict}, we get
$$|s_N-t_0|^2+\|Z_N-X_0\|_2^2\le(\hat U-\Phi)(t_0,X_0)-(\hat U-\Phi)(t_0,X_{\bm x_N})+2e_N\xrightarrow{N\to+\infty}0 .$$
Hence, for $N$ large, $(s_N,\bm z_N)$ is an interior point of $K_N$ with $s_N>0$. The subsolution property of $V^N$ in light of \eqref{Hdisc} then yields
\be\label{discineq}
\partial_t\Phi(s_N,Z_N)+\hat H\Big(Z_N,\,\mathbb E\big[D\Phi(s_N,Z_N)\,|\,\mathcal F_N\big]\Big)\le0 .
\ee
Since $(s_N,Z_N)\to(t_0,X_0)$, we have $D\Phi(s_N,Z_N)\to Y_0$. Moreover, since $Y_N$ is $\mathcal F_N$-measurable and the conditional expectation is the best $L^2$-approximation by $\mathcal F_N$-measurable functions we get
\begin{align*}
\big\|\mathbb E[D\Phi(s_N,Z_N)|\mathcal F_N]-Y_0\big\|_2 &\le \big\|D\Phi(s_N,Z_N)-Y_0\big\|_2+\big\|\mathbb E[D\Phi(s_N,Z_N)|\mathcal F_N]-D\Phi(s_N,Z_N)\big\|_2\\
&\le 2\|D\Phi(s_N,Z_N)-Y_0\|_2+\|Y_N-Y_0\|_2\xrightarrow{N\to+\infty}0 .
\end{align*}
Letting $N\to+\infty$ in \eqref{discineq} and using Assumption (A1)(1), we obtain the subsolution inequality
$$\partial_t\Phi(t_0,X_0)+\hat H(X_0,D\Phi(t_0,X_0))\le0.\qedhere$$
\end{proof}

\section{Example Hamiltonians}\label{ex1}

\nd
In this section of the appendix we present the proofs of the examples provided in Remark \ref{asrem}.

\subsection{Convex Hamiltonians}
\begin{lemma}\label{lem:SCA}
    Assume that $H$ satisfies (A1) and is such that $p\mapsto H(x,p,m)-\frac{\nu}{2}|p|^2$ is convex, for some $\nu>0$. For every $\eta>0$, there exists a constant $\psi_\eta>0$, depending only on $\eta,\nu,c_H,C_H, K_0,k_0$, such that the following holds for all $\mu\in(0,1]$ and $\theta\in[0,\eta/2]$.
    \begin{enumerate}
        \item[(i)] If $x\in\T^d$, $m\in\mathcal{P}_2$ and $p,q\in\R^d$ satisfy $H(x,p+q,m)\ge \eta$ and $H(x,q,m)\le \theta$, then
        \be\label{pSCA}
        \mu H(x,p/\mu+q,m)\ge H(x,p+q,m)+(\psi_\eta-\theta)(1-\mu).
        \ee
        \item[(ii)] If $X,P,Q\in\mathcal{H}$ satisfy $\hat{H}(X,P+Q)\ge \eta$ and $\hat{H}(X,Q)\le \theta$, then
        \be \label{iSCA}\tag{\text{i-SCA}}
        \mu \hat{H}(X,P/\mu+Q)\ge \hat{H}(X,P+Q)+(\psi_\eta-\theta)(1-\mu).
        \ee
    \end{enumerate}
    In particular, since $H+\lambda$ is also $\nu$-strongly convex in $p$ and satisfies (A1) with $k_0$ and $K_0$ replaced by $k_0-\lambda\ge0$ \footnote{Recall Proposition \ref{ergodicconstant}.} and $K_0+\lambda$, applying (ii) to $H+\lambda$ shows that $H$ satisfies \eqref{SCA}.
    % In particular, for $\theta=0$ we recover the (SCA) condition: if $H(x,p+q,m)\ge\eta$ and $H(x,q,m)\le 0$, then $\mu H(x,p/\mu+q,m)\ge H(x,p+q,m)+\psi_\eta(1-\mu)$ for all $\mu\in(0,1]$, and similarly for $\hat{H}$.
\end{lemma}

\begin{proof}
    We prove (ii); the proof of (i) is identical, with expectations removed and the $L^2$-norm $\|\cdot\|_2$ replaced by $|\cdot|$. There is nothing to prove if $\mu=1$, so let $\mu\in (0,1)$. Since $p\mapsto H(x,p,m)-\frac{\nu}{2}|p|^2$ is convex, for all $a,b\in\R^d$ we have
    $$H(x,\mu a+(1-\mu)b,m)\le \mu H(x,a,m)+(1-\mu)H(x,b,m)-\frac{\nu}{2}\mu(1-\mu)|a-b|^2.$$
    We apply this for every $\omega\in\Omega$ with $x=X(\omega)$, $m=\mathcal{L}(X)$, $a=P(\omega)/\mu+Q(\omega)$ and $b=Q(\omega)$, so that $\mu a+(1-\mu)b=P(\omega)+Q(\omega)$ and $a-b=P(\omega)/\mu$. This gives
    $$\mu H\left(X,\frac{P}{\mu}+Q,\mathcal{L}(X)\right)\ge H(X,P+Q,\mathcal{L}(X))-(1-\mu)H(X,Q,\mathcal{L}(X))+\frac{\nu}{2}\frac{1-\mu}{\mu}|P|^2.$$
    All terms are integrable by (A1)(2). Taking expectations and using $\hat H(X,Q)\le \theta$ and $\mu\le 1$, we get
    \be\label{SCAint}
    \mu \hat{H}\left(X,\frac{P}{\mu}+Q\right)\ge \hat{H}(X,P+Q)-(1-\mu)\theta+\frac{\nu}{2}(1-\mu)\|P\|_2^2.
    \ee
    It remains to bound $\|P\|_2$ from below by a constant depending only on $\eta$. By Assumption (A1)(2), $c_H\|Q\|_2^2-k_0\le \hat{H}(X,Q)\le \theta\le \eta/2$, hence $\|Q\|_2\le R_\eta:=(\frac{k_0+\eta/2}{c_H})^{1/2}$. By Assumption (A1)(1) (with $x=x'$, $m=m'$) and the Cauchy--Schwarz inequality, we have
    $$\frac{\eta}{2}\le \eta-\theta\le \hat{H}(X,P+Q)-\hat{H}(X,Q)\le C_H\,\mathbb{E}\big[(1+|P+Q|+|Q|)|P|\big]\le C_H\left(1+\|P\|_2+2R_\eta\right)\|P\|_2.$$
    Let $r_\eta>0$ be the unique positive root of $C_H(1+r+2R_\eta)r=\eta/2$. Since the right-hand side is increasing in $\|P\|_2$, we get $\|P\|_2\ge r_\eta$. The conclusion follows from \eqref{SCAint} with $\psi_\eta:=\frac{\nu}{2}r_\eta^2$.
\end{proof}

\begin{remark}
    The pointwise (SCA) condition alone does not imply the integral version (ii), since the conditions $\hat{H}(X,P+Q)\ge\eta$ and $\hat{H}(X,Q)\le 0$ do not imply $H(X,P+Q,\mathcal{L}(X))\ge\eta$ and $H(X,Q,\mathcal{L}(X))\le 0$ almost surely. 
    %This is why we derive both versions directly from the uniform convexity of $H$. The key point is that the uniform convexity inequality holds pointwise without any preconditions. The preconditions are only used to bound $\|P\|_2$ from below.
\end{remark}

\subsection{Non-convex Hamiltonians}
 
\begin{lemma}\label{lem:gameISCA}
Let $H$ be the  Hamiltonian given in Remark \ref{asrem}(ii), that is
$$H(x,p,m)=\frac{1}{2}|p|^2-\rho \chi(m)(|p|-s_0)_+-F(m),$$
where $\rho,s_0>0$ with $\e:=\frac{\rho}{s_0}<\frac{1}{6}$,
$$l(m):=\int_{\T^d}f(x) m(dx),\;\; F(m):= (l(m)-1)_+,\;\;\text{and }\;\;\chi(m)=\min\{ 1, (1-l(m))_+\},\;\;m\in\mathcal{P}_2,$$
for some function $f\in C^1(\T^d)$ so that $\min f<0<1<\max f$. Then,
\begin{enumerate}
\item[(i)] $H$ satisfies (A1);
\item[(ii)] the ergodic constant of Proposition~\ref{ergodicconstant} is $\lambda=0$;
\item[(iii)] there exists $(x,m)\in \T^d\times \mathcal{P}_2$ so that the map $p\mapsto H(x,p,m)$ is not convex;
\item[(iv)] $H$ satisfies \eqref{SCA} with $c_*:=\frac{1}{1-2\varepsilon}\ge 1$.
\end{enumerate}
\end{lemma}
 
\begin{proof}
Let $\phi(p):=-\rho\,(|p|-s_0)_+$, so that
$H(x,p,m)=\frac12|p|^2+\chi(m)\phi(p)-F(m)$. It is easy to check that $\phi$ is concave and $\rho$-Lipschitz, and $\phi\le0$. In addition, $|\phi(p)-\phi(p')|\le\varepsilon\max\{|p|,|p'|\}\,|p-p'|$, for all $p,p'\in \R^d$. It is also straightforward to check that for $\mu\in(0,1]$
\be\label{P3}
\big|\mu\,\phi(q+p/\mu)-\phi(q+p)+(1-\mu)\phi(q)\big|
\le\varepsilon(1-\mu)\big(3|p|^2/\mu+|q|^2\big),\;\;\text{for all }p,q\in\mathbb{R}^d.
\ee
Since $f\in C^1(\T^d)$, the functions $l,F, \chi$ are ${\bf d}_2$-Lipschitz with constant the Lipschitz constant of $f$.
\vspace{2mm}

\nd
\textit{Proof of (i).} The boundedness of $f$ implies 
\be\label{game:growth}
\Big(\frac12-\varepsilon\Big)|p|^2-F(m)\le H(x,p,m)\le\frac12|p|^2,
\ee
and since $F$ is also bounded, (A1)(2) follows. This also implies that
\begin{equation}\label{game:min}
\min_{p\in\mathbb{R}^d}H(x,p,m)=H(x,0,m)=-F(m).
\end{equation}
(A1)(1) follows from the properties of $\chi,F,l,\phi$ discussed at the beginning.
\vspace{2mm}

\nd
\textit{Proof of (ii).} Let $V$ be a ${\bf d}_2$-Lipschitz L-viscosity solution of
\eqref{EHJB} with constant $\lambda$, given by
Proposition~\ref{ergodicconstant}. Since $\mathcal{P}_2$ is compact, $V$ attains its minimum at
some $m_1$, and every $X_1\in\mathcal{H}$ with $\mathcal{L}(X_1)=m_1$ is a
minimum point of $\widehat V$. Testing the supersolution property
(Definition~\ref{evisc.}) with a constant test function gives
$\lambda+\widehat H(X_1,0)\ge0$, that is, $\lambda\ge F(m_1)\ge0$.
\vspace{1mm}

\nd
We show the opposite inequality $\lambda\le0$. Recall from the proof of Proposition~\ref{ergodicconstant}
that $\lambda$ is the limit, along a subsequence, of the ergodic constants
$\lambda^N$ of the $N$-particle problems
\[
\lambda^N+\frac1N\sum_{i=1}^NH\big(x^i,ND_{x^i}V^N(\bm x),m^N_{\bm x}\big)=0,\quad\text{in }(\mathbb{T}^d)^N,
\]
which admit Lipschitz viscosity solutions $V^N$. Since $V^N$ is Lipschitz, it is differentiable a.e. At points of differentiability, the
subsolution property gives 
$$\lambda^N+\frac1N\sum_{i=1}^NH\big(x^i,ND_{x^i}V^N(\bm x),m^N_{\bm x}\big)\le 0\implies \lambda^N+\min_{P\in (\R^d)^N}\frac1N\sum_{i=1}^NH\big(x^i,NP^i,m^N_{\bm x}\big)\le 0,$$
for a.e.\ $\bm x\in (\T^d)^N$. By \eqref{game:min}, this implies $\lambda^N-F(m^N_{\bm x})\le 0$, for a.e.\ $\bm x\in (\T^d)^N$, and therefore for all $\bm x\in (\T^d)^N$ by
continuity. 
Choose $x_*\in\mathbb{T}^d$ with $f(x_*)=\min f$ and let $\bm x=(x_*,\dots,x_*)$. Then
$m^N_{\bm x}=\delta_{x_*}$ and $l(\delta_{x_*})=\min f<0$, so $F(\delta_{x_*})=0$.
Hence $\lambda^N\le0$, and so $\lambda\le0$.
\vspace{2mm}

\nd
\textit{Proof of (iii).} Let $e_1=(1,0,\ldots,0)\in \R^d$ and $(x,m)\in \T^d\times\mathcal{P}_2$ so that $\chi(m)>0$ and hence $F(m)=0$. The function
$k(s):=H(x,se_1,m)=\frac12s^2-\rho\chi(m)(s-s_0)_+$, $s\ge0$, satisfies
\[
\lim_{s\to s_0^-}k'(s)=s_0>s_0-\rho\chi(m)=\lim_{s\to s_0^+}k'(s).
\]
A convex function cannot have a downward jump in its derivative, so
$p\mapsto H(x,p,m)$ is not convex.

\nd
\textit{Proof of (iv).} By (ii), $\lambda=0$. Fix $\eta>0$, $\mu\in(0,1)$ (the
case $\mu=1$ is trivial), $\theta\in[0,\eta/2]$, and $X,P,Q\in\mathcal{H}$
with $\widehat H(X,P+Q)\ge\eta$ and $\widehat H(X,Q)\le\theta$. Let
$m:=\mathcal{L}(X)$ and set 
$$k(p):=\frac12|p|^2+\chi(m)\phi(p)-F(m).$$
Since $H$ does not depend on $x$, we have $\widehat H(X,R)=\mathbb{E}[k(R)]$ for every
$R\in\mathcal{H}$.

\nd
\emph{(iv-a) Lower bound on $\|P\|_2$.} By \eqref{game:growth} and the boundedness of $F$,
we have $(\frac12-\varepsilon)\|Q\|_2^2-M\le\theta\le\eta/2$, for some constant $M>0$; hence $\|Q\|_2\le R_\eta:=\sqrt{\frac{\eta/2+M}{1/2-\e}}$. By the properties of $\phi, \chi, F$, we have
$|k(p)-k(p')|\le(|p|+|p'|)|p-p'|$, so
\[
\frac\eta2\le\eta-\theta\le\mathbb{E}\big[k(P+Q)-k(Q)\big]
\le\mathbb{E}\big[(|P|+2|Q|)|P|\big]\le\|P\|_2\big(\|P\|_2+2R_\eta\big).
\]
Letting $r_\eta$ be the positive root of $r(r+2R_\eta)=\eta/2$, this gives
$\|P\|_2\ge r_\eta$.

\nd
\emph{(iv-b) Conclusion} For $p,q\in\mathbb{R}^d$ set
$a:=q+p/\mu$ and $r:=q+p=\mu a+(1-\mu)q$. We use the identity
\[
\mu|a|^2+(1-\mu)|q|^2-|r|^2=\mu(1-\mu)|a-q|^2=\frac{(1-\mu)|p|^2}{\mu},
\]
the cancellation of the constant $-F(m)$, and \eqref{P3} to get
\begin{align*}
\mu k(a)-k(r)+(1-\mu)k(q)
&=\frac{(1-\mu)|p|^2}{2\mu}+\chi(m)\big(\mu\phi(a)-\phi(r)+(1-\mu)\phi(q)\big)\\
&\ge(1-\mu)\Big(\Big(\frac12-3\varepsilon\chi(m)\Big)\frac{|p|^2}{\mu}
-\varepsilon\chi(m)|q|^2\Big).
\end{align*}
Applying this inequality for $p=P(\omega)$ and $q=Q(\omega)$,
and taking expectations. Using $\widehat H(X,Q)\le\theta$, $\mu\le1$ and
$\frac12-3\varepsilon\chi(m)\ge\frac12-3\varepsilon>0$, we get
$$\mu\,\widehat H\Big(X,\frac P\mu+Q\Big)\ge\widehat H(X,P+Q)
+(1-\mu)\Big(\Big(\frac12-3\varepsilon\Big)\|P\|_2^2
-\big(\theta+\varepsilon\chi(m)\|Q\|_2^2\big)\Big).$$
We claim that $\theta+\varepsilon\chi(m)\|Q\|_2^2\le c_*\theta$. Indeed, if $\chi(m)=0$, this is clear, since $c_*\ge1$.  If $\chi(m)>0$, then $F(m)=0$ and the left side of \eqref{game:growth} gives
$(\frac12-\varepsilon)\|Q\|_2^2\le\theta$. Hence
$\varepsilon\chi(m)\|Q\|_2^2\le\frac{2\varepsilon}{1-2\varepsilon}\theta$, and
the claim follows.
\vspace{1mm}

\nd
Combining the claim with (iv-a) yields
\[
\mu\,\widehat H\Big(X,\frac P\mu+Q\Big)\ge\widehat H(X,P+Q)
+\big(\psi_\eta-c_*\theta\big)(1-\mu),\;\;\text{where }\;\; \psi_\eta:=(\tfrac12-3\varepsilon)r_\eta^2,
\]
which, since $\lambda=0$, is \eqref{SCA}.
\end{proof}

\bibliographystyle{alpha-author}
\bibliography{references}

\end{document}